\input amstex
\documentstyle{amsppt}
\magnification=\magstep1
\baselineskip 20pt
\topmatter
\title 
On the complete classification of extremal log 
Enriques surfaces
\endtitle
\author
Keiji Oguiso and De-Qi Zhang
\endauthor 
\abstract
We show that there are exactly, up to isomorphisms, seven
rational extremal log Enriques surfaces  Z  and construct all of them;
among them types  $D_{19}$  and  $A_{19}$  have been shown
of certain uniqueness by M. Reid.  We also prove that the
(degree 3 or 2) canonical covering of each of these seven  Z  
has either  $X_3$  or  $X_4$  as its minimal resolution.  Here
$X_3$ (resp. $X_4$)  is the unique K3 surface with
Picard number 20 and discriminant 3 (resp. 4), which
are called the most algebraic K3 surfaces by Vinberg
and have infinite automorphism groups
(by Shioda-Inose and Vinberg).
\endabstract
\endtopmatter
\document 
\head
Introduction
\endhead

\par 
Throughout this paper, we work over the complex number field $\Bbb C$. 
A normal projective surface $Z$ with at worst 
quotient singularities is called {\it a log Enriques surface} 
if the canonical Weil divisor  $K_{Z}$  is numerically equivalent
to zero and if the irregularity  dim $H^{1}(Z, \Cal O_{Z}) = 0$ 
[Z1, (1.1)].  Note that a log Enriques surface is irrational 
if and only if it is a K3 or Enriques surface with at worst 
Du Val singular points, and also we can regard log Enriques 
surfaces as degenerations or generalizations of K3 and Enriques surfaces. 

\par 
Rational log Enriques surfaces also appear as base spaces  $W$  of
elliptically fibred Calabi-Yau threefolds
$\Phi_{|D|} : X \rightarrow W$  with  $D. c_2(X) = 0$ [O1].
On the other hand, a special case of [Al, Theorem 3.9] says that
there are only finitely many deformation types of minimal resolutions
of rational log Enriques surfaces.  This also shows the sort of
feasibility to classify them all.

\par
Since the minimal partial 
resolution of the Du Val singular points of a log Enriques surface is 
again a log Enriques surface of the same canonical index (see below for
the definition), we assume throughout this paper that 
a log Enriques surface has no Du Val singular points.  

\par
Let $Z$ be a log Enriques surface and 
$$ 
I := I(Z) = \text{min}\{n \in \Bbb Z_{> 0} \vert \Cal O_{Z}(nK_{Z}) 
\simeq \Cal O_{Z}\}
$$ 
{\it the canonical index} of $Z$. {\it The canonical covering} 
of $Z$ is then defined as
$$ 
\pi : \overline{S} := Spec_{\Cal O_{Z}} 
(\oplus_{i=0}^{I-1} \Cal O_{Z}(-iK_{Z})) \rightarrow Z.
$$ 
It follows from [K, Z1] that 
\roster 
\item 
$\overline{S}$ is either a projective K3 surface with at worst 
Du Val singularities or an abelian surface, 
\item 
$\pi$ is a finite, cyclic Galois cover of degree $I$ and is \'etale 
over $Z - \text{Sing}Z$, and 
\item 
the Galois group $\text{Gal} (\overline{S}/Z) \simeq \Bbb Z/I \Bbb Z$ 
acts faithfully on $H^{0}(\Cal O_{\overline{S}}(K_{\overline{S}}))$, 
that is, there exists a generator $g$ of $\text{Gal}(\overline{S}/Z)$ 
such that 
$g^{*}\omega_{\overline{S}} = \zeta_{I}\omega_{\overline{S}}$, 
where $\zeta_{I} = \text{exp}(2\pi i/I)$ and 
$\omega_{\overline{S}}$ is a nowhere vanishing regular 2-form on 
$\overline{S}$.
\endroster

\par
One interesting problem is to determine all possible
canonical indices; in this aspect, [Bl] has shown that 
the canonical index is always less than or equal to 21 (see also [Z1,2]).
On the other hand, in [Z1,2] for each prime number  
$p \le 19$, we have constructed 
a rational log Enriques surface  $Z_p$  of index  $p$, with
the canonical covering  $\pi : Y_p \rightarrow Z_p$, the Galois
group  $G = Gal(Y_p/Z_p)$  and the minimal resolution
$X_p \rightarrow Y_p$, while
in [OZ3] we have shown that for each  $p = 13, 17, 19$  the 
pair  $(X_p, G)$  is unique up to isomorphisms.

\par
Let $\nu : S \rightarrow \overline{S}$ be the minimal resolution 
of $\overline{S}$ and $\Delta_{Z}$ the exceptional divisor of $\nu$. 
Then $\Delta_{Z}$ is a disconnected sum of divisors of Dynkin's type, 
$(\oplus A_{\alpha}) \oplus (\oplus D_{\beta}) \oplus (\oplus E_{\gamma})$. 
Then, by abuse of language, we say that a log Enriques 
surface $Z$ or the exceptional divisor $\Delta_{Z}$ is of type 
$(\oplus A_{\alpha}) \oplus (\oplus D_{\beta}) \oplus (\oplus E_{\gamma})$. 
We define $\text{rank}\Delta_{Z}$ as the rank of the sublattice of the 
N\'eron Severi lattice $NS(S) = \text{Pic} Z$ generated by the irreducible components of $\Delta_{Z}$. 
Note that  $\text{rank}\Delta_{Z}$ is the number 
$\sum \alpha + \sum \beta + \sum \gamma$  of 
the exceptional curves and satisfies  
$$
\text{rank}\Delta_{Z} \leq \text{rank}NS(S) -1 \leq 19. 
$$

\par
Our standpoint here is, as in previous [OZ1, 2], 
to consider  $\text{rank}\Delta_{Z}$ as an invariant 
measuring how ``big'' $\text{Sing}(Z)$  is.

\par
\proclaim{Definition} A rational log Enriques surface $Z$ is said to be 
extremal if $\text{rank} \Delta_{Z} = 19$, the maximal possible value.
\endproclaim 

\par
Note that the minimal resolution $S$ of the canonical cover of an 
extremal log Enriques surface is necessarily a singular K3 surface, 
that is, a smooth K3 surface with maximal possible Picard number 20,
in the terminology of [SI].  Thus, it is very natural to
ask whether we can show the uniqueness of each extremal type,
up to isomorphisms [see Question 2 below].

\par
In [OZ2], we have determined isomorphism classes of rational
log Enriques surfaces of type
$D_{18}$  (one class only) or  
$A_{18}$  (two classes), while in [OZ1]  
we gave an affirmative answer to the following question raised to 
the second author by I. Naruki and M. Reid when they saw
the examples of rational log Enriques surfaces of type
$D_{19}$  or  $A_{19}$  in [Z1].

\par
\proclaim{Question 1} Are the rational log Enriques surfaces 
of type $D_{19}$ and of type $A_{19}$ unique respectively
up to isomorphisms?
\endproclaim 

\par
This question is now naturally generalised to the following:
\proclaim{Question 2} How about the extremal log Enriques surfaces?
\endproclaim

\par
The main purpose of this paper is to give a complete 
answer to Question 2:

\par 
\proclaim{Main Theorem} 

\par
(1)  [Restriction] Each extremal log Enriques surface has 
one of the following seven types:
\par \hskip 2pc 
$D_{19}$, $D_{16} \oplus A_{3}$, $D_{13} \oplus A_{6}$, 
$D_{7} \oplus A_{12}$, $D_{7} \oplus D_{12}$, 
$D_{4} \oplus A_{15}$, or $A_{19}$. 

\par
(2)  [Existence] Conversely, for each type $\Delta$ given in (1),
there exists an extremal log Enriques surface of type $\Delta$.
(See Prop. 1.7.)

\par
(3)  [Uniqueness] Extremal log Enriques surfaces are isomorphic if and only if 
their types are the same.

\par
In particular, there exist exactly seven extremal log Enriques surfaces up 
to isomorphisms.
\endproclaim 

\par
In section 1, we explicitly construct an extremal log Enriques surface 
of each type given in (1) via Shioda-Inose's pairs $(S_{3}, <g_{3}>)$ and 
$(S_{2}, <g_{2}>)$, that is, pairs of the singular K3 surfaces 
with two smallest discriminants 3 and 4 and their certain automorphism group 
of order 3 and 2 respectively.  (See \S 1 for the detail.)
As in [OZ1], the basic strategy of the proof here for the Main Theorem is 
to reduce problems of an extremal log Enriques surface to those of a singular 
K3 surfaces via the canonical covering and its minimal resolution,
the so-called Godeaux approach.
In section 2, we show the following proposition, which determines 
extremal log Enriques surfaces except for an ambiguity of 
the exceptional divisor $\Delta_{Z}$ of $S \rightarrow \overline{S}$,
and is one of the cores of this paper:
 
\par
\proclaim{Proposition (cf. Proposition (2.2))} 
Let $Z$ be an extremal log Enriques surface, $\overline{S} \rightarrow Z$ 
the canonical cover of $Z$ and $S \rightarrow \overline{S}$ 
the minimal resolution of $\overline{S}$. Let $<g>$ be the 
automorphism group of $S$  
induced by the Galois group of $\overline{S} \rightarrow Z$. 
Then, the pair $(S, <g>)$ is isomorphic to either one of 
Shioda-Inose's pairs $(S_{3}, <g_{3}>)$ or $(S_{2}, <g_{2}>)$. 
In particular, the canonical index of an extremal log Enriques surface is 
either 3 or 2.
\endproclaim  

\par
The hardest part of this proposition is the determination 
of the canonical indices of extremal log Enriques surfaces.
For this, we need some detailed analysis of the fixed
locus  $S^{<g>}$  based on Atiyah-Singer-Segal's 
holomorphic Lefschetz fixed point formula [AS1,2] and 
the usual topological 
Lefschetz fixed point formula (see eg. [U]).
This analysis, which describes the fixed locus
of an order 6  automorphism  $\tau$  on a K3 surface  $T$
with  $\tau^* \omega = \zeta_6 \omega_T$,
will be applicable to quite general cases.
After proving  $I = 2$  or  $3$, we apply 
the characterisation of Shioda-Inose's pairs 
$(S_{3}, <g_{3}>)$ and $(S_{2}, <g_{2}>)$
(Theorems (1.3) and (1.6)) to conclude 
Proposition(2.2).

\par
In section 3, studying $\Delta$ as a sublattice of $NS(S)$, 
we show the Main Theorem (1). In section 4, we prove the Main 
Theorem (3) along the strategy given in (4.1).  

\par \vskip 5pt
{\bf Acknowledgement.}  The present version of this paper
was completed during the first author's visit to Singapore,
and the authors would like to express their gratitude to
the JSPS programme and the National University of Singapore
for the financial support.  The authors also would like 
to thank the referee for suggestions which make the paper
more comprehensive.

\par
\head 
Notation 
\endhead 
For an automorphism group $G$  and its element  $g$  
of a smooth surface or a curve $X$, we set
\par \noindent
$X^{g} := \{x \in X \vert g(x) = x\}$, the fixed locus 
of an element $g$ of $G$,
\par \noindent 
$X^{[G]} := \cup_{g \in G - \{id\}} X^{g} = 
\{ x \in X \vert g(x) = x$ for some $g \in G - \{id\}\}$. 
Note that  $X^g$  is a smooth algebraic set.

\par \noindent 
A curve $C$  on a surface  $S$ is said to be 
$g-$stable (resp. $g-$fixed) if 
$g(C) = C$ (resp. $C \subset S^{g}$). We call $P \in S^{g}$ an 
isolated point if $P$ is not contained in any $g-$fixed curves. 
\par \noindent  
We denote by $\omega_{S}$ a nowhere vanishing holomorphic 2-form of a 
K3 surface $S$.
\par \noindent
We denote by $\zeta_{I} = \text{exp}(2\pi i/I)$,  
a specified primitive $I-$th root of unity. 

\par
Let $Z$ be an extremal log Enriques surface. We set:
\par \noindent
$I = I(Z)$ the canonical index of $Z$;
\par \noindent
$\pi : \overline{S} \rightarrow Z$ the canonical cover of $Z$;
\par \noindent
$g$ the generator of $Gal(\overline{S}/Z) (\simeq \Bbb Z/I \Bbb Z)$ 
such that
$g^{*}\omega_{\overline{S}} = \zeta_{I}\omega_{\overline{S}}$;
\par \noindent
$\nu : S \rightarrow \overline{S}$ the minimal resolution of
$\overline{S}$;
\par \noindent
$\Delta = \Delta_{Z}$ the exceptional divisor of $\nu$;
\par \noindent
$\Delta = \Delta_{1} \coprod \cdots \coprod \Delta_{r}$
the decomposition of $\Delta$ into the connected components.
\par \noindent
For the simplicity of notation, we also denote by the same letter
$g$ the induced action of $g$ on $S$.

\head 
\S 1. Construction of extremal log Enriques surfaces from Shioda-Inose's pairs 
\endhead
First, we recall definitions and some properties of Shioda-Inose's pairs $(S_{3}, <g_{3}>)$ and $(S_{2}, <g_{2}>)$  from [SI]
and [OZ1].  These pairs will play essential roles throughout 
this paper.  Next we construct extremal log Enriques surfaces 
of all types shown in the Main Theorem (1). 
This will complete main Theorem (2), the existence part. 

\par 
\definition{Definition (1.1) [OZ1, Example 1]}  
Let $E_{\zeta_{3}}$ be the elliptic curve of period $\zeta_{3}$.  
Let $\overline{S_{3}} := E_{\zeta_{3}}^{2}/<\text{diag}(\zeta_{3}, \zeta_{3}^{2})>$ be the quotient surface of the product $E_{\zeta_{3}}^{2}$ by $<\text{diag}(\zeta_{3}, \zeta_{3}^{2})>$ and  $S_{3} \rightarrow \overline{S_{3}}$ the minimal resolution of $\overline{S_{3}}$.
\par 
Let $g_{3}$ be the automorphism of $S_{3}$  
(of order 3) induced by the action $\text{diag}(\zeta_{3}, 1)$ on $E_{\zeta_{3}}^{2}$. 
We call the pair  $(S_{3}, g_{3})$  or  $(S_{3}, <g_{3}>)$
Shioda-Inose's pair of discriminant 3.

\enddefinition
\remark{Remark}
It is shown in [SI] that $S_{3}$ is uniquely determined by the property that $S_{3}$ is a singular K3 surface whose transcendental lattice is of 
discriminant 3. This is the source of the name given in (1.1).
\endremark

\par \vskip 5pt
The next Proposition and Theorem are shown in [SI, OZ1].
\proclaim{Proposition (1.2)} 
Let $(S_{3}, g_{3})$ be the Shioda-Inose's pair of discriminant 3. Then,
\roster 
\item 
$S_{3}$ contains 24 rational curves: 
$F_{1}$, $F_{2}$, $F_{3}$ coming from 
$(E_{\zeta_{3}})^{\zeta_{3}} \times E_{\zeta_{3}}$, 
$G_{1}$, $G_{2}$, $G_{3}$ coming from 
$E_{\zeta_{3}} \times (E_{\zeta_{3}})^{\zeta_{3}}$, and 
$E_{ij}$, $E_{ij}'$ ($1 \leq i,j \leq 3$) the exceptional curves 
arising from the nine Du Val singular points (of Dynkin type 
$A_{2}$) of $\overline{S_{3}}$ 
(see [SI, Fig.2] or Figure 1 at the end of this paper 
for the configuration of these 24 curves), 
\item 
$g_{3}^{*}\omega_{S_{3}} = \zeta_{3}\omega_{S_{3}}$ and each of the 24 curves in (1) 
is $g_{3}-$stable, 
\item
$S_{3}^{g_{3}} = (\coprod_{i=1}^{3} F_{i}) \coprod (\coprod_{j=1}^{3} G_{j}) 
\coprod (\coprod_{i,j = 1}^{3} \{P_{ij}\})$, where $P_{ij}$ denote the point 
$E_{ij} \cap E_{ij}'$, and 
\item
$g_{3} \circ \varphi = \varphi \circ g_{3}$ for all $\varphi \in \text{Aut}(S_{3})$.
\endroster

\par
\endproclaim 
\proclaim{Theorem (1.3)} 
Let $(S, g)$ be a pair of a smooth K3 surface and 
an automorphism $g$ of $S$. 
Assume that $(S, g)$ satisfies the following four conditions:
\roster
\item 
$g^{3} = id$,
\item 
$g^{*}\omega_{S} = \zeta_{3}\omega_{S}$,
\item 
$S^{g}$ consists of only rational 
curves and isolated points, and 
\item 
$S^{g}$ contains at least six rational curves.
\endroster
Then $(S, g) = (S_{3}, g_{3})$ up to equivariant isomorphisms. Moreover, 
$S^{g}$ consists of exactly six rational curves and nine isolated points.
\endproclaim 

\par 
\definition{Definition (1.4) [OZ1, Example 2]}  
Let $E_{\zeta_{4}}$ be the elliptic curve of period $\zeta_{4}$.  
Let $\overline{S_{2}} := E_{\zeta_{4}}^{2}/<\text{diag}(\zeta_{4}, \zeta_{4}^{3})>$ be the quotient surface of the product $E_{\zeta_{4}}^{2}$ by $<\text{diag}(\zeta_{4}, \zeta_{4}^{3})>$ and  $S_{2} \rightarrow \overline{S_{2}}$ the minimal resolution of $\overline{S_{2}}$.
\par 
Let $g_{2}$ be the involution of $S_{2}$  
induced by the action $\text{diag}(-1, 1)$ on $E_{\zeta_{4}}^{2}$. 
We call the pair $(S_{2}, g_{2})$ Shioda-Inose's pair of discriminant 4.
\enddefinition
\remark{Remark}
It is also shown in [SI] that $S_{2}$ is uniquely determined by the property that $S_{2}$ is a singular K3 surface whose transcendental lattice is of 
discriminant 4. 
\endremark

\par \vskip 5pt
Proposition(1.5) and Theorem(1.6) below are also shown in [SI, OZ1].

\par
\proclaim{Proposition (1.5)} 
Let $(S_{2}, g_{2})$ be Shioda-Inose's pair of discriminant 4. Then,
\roster 
\item 
$S_{2}$ contains 24 rational curves: 
$F_{1}$, $F_{2}$, $F_{3}$ coming from 
$(E_{\zeta_{4}})^{[<\zeta_{4}>]} \times E_{\zeta_{4}}$, 
$G_{1}$, $G_{2}$, $G_{3}$ coming from 
$E_{\zeta_{4}} \times (E_{\zeta_{4}})^{[<\zeta_{4}>]}$,
$E_{11}' + H_{11} + E_{11}$, $E_{13}' + H_{13} + E_{13}$, 
$E_{31}' + H_{31} + E_{31}$, $E_{33}' + H_{33} + E_{33}$, the 
exceptional curves arising from the four Du Val singular points 
(of Dynkin type $A_{3}$) of $\overline{S_{2}}$ , and
$E_{12}$, $E_{22}$, $E_{32}$ $E_{21}'$, $E_{22}'$, $E_{23}'$, 
the exceptional curves 
arising from the six Du Val singular points (of type 
$A_{1}$) of $\overline{S_{2}}$ 
(see [SI, Fig.3] or Figure 2 at the end of this paper 
for the configuration of these 24 curves), 
\item 
$g_{2}^{*}\omega_{S_{2}} = - \omega_{S_{2}}$ and each of the 24 curves in (1) 
is $g_{2}-$stable, 
\item
$S_{2}^{g_{2}} = (\coprod_{i=1}^{3} F_{i}) \coprod 
(\coprod_{j=1}^{3} G_{j}) 
\coprod (\coprod_{i,j = 1}^{3} H_{ij}$), and
\item 
$g_{2} \circ \varphi = \varphi \circ g_{2}$ for all $\varphi \in \text{Aut}(S_{2})$. 
\endroster
\endproclaim 
\proclaim{Theorem (1.6)} 
Let $(S, g)$ be a pair of a smooth K3 surface and 
an automorphism $g$ of $S$. 
Assume that $(S, g)$ satisfies the following four conditions:
\roster
\item 
$g^{2} = id.$,
\item 
$g^{*}\omega_{S} = -\omega_{S}$,
\item 
$S^{g}$ consists of only rational curves, and 
\item 
$S^{g}$ contains at least ten rational curves.
\endroster
Then $(S, g) = (S_{2}, g_{2})$ up to equivariant isomorphisms. Moreover, 
$S^{g}$ consists of exactly ten rational curves.
\endproclaim  
Now, using the notation in (1.1), (1.2), (1.4), (1.5) and tracing out
Figures 1 and 2, we can easily construct an extremal 
log Enriques surface of 
each type given in the Main Theorem (1) as follows:

\par
\proclaim{Proposition (1.7)} 

\roster 
\item 
Let $\Delta (i)$ ($1 \leq i \leq 6$) be the divisors on $S_{3}$ 
defined by:

\par \vskip 5pt \noindent
$\Delta (1) = E_{11} + E_{21} + G_{1} + E_{31} + E_{31}' + F_{3} + 
E_{33}' + E_{33} + G_{3} + E_{23} + E_{23}' + F_{2} + 
E_{22}' + E_{22} + G_{2} + E_{12} + E_{12}' + F_{1} + E_{13}'$ 
(of Dynkin type $D_{19}$);

\par \vskip 5pt \noindent
$\Delta (2) = (E_{11}' + E_{12}' + F_{1} + E_{13}' + E_{13} + G_{3} + 
E_{23} + E_{23}' + F_{2} + E_{22}' + E_{22} + G_{2} + 
E_{32} + E_{32}' + F_{3} + E_{33}') + (E_{21} + G_{1} + E_{31})$ 
(of Dynkin type $D_{16} \oplus A_{3}$);

\par \vskip 5pt \noindent
$\Delta (3) = (E_{12}' + E_{13}' + F_{1} + E_{11}' + E_{11} + G_{1} + 
E_{21} + E_{21}' + F_{2} + E_{22}' + E_{22} + G_{2} + 
E_{32}) + (E_{31}' + F_{3} + E_{33}' + E_{33} + G_{3} + E_{23})$ 
(of Dynkin type $D_{13} \oplus A_{6}$);

\par \vskip 5pt \noindent
$\Delta (4) = (E_{11}' + E_{12}' + F_{1} + E_{13}' + E_{13} + G_{3} + 
E_{23}) + (E_{33}' + F_{3} + E_{32}' + E_{32} + G_{2} + 
E_{22} + E_{22}' + F_{2} + E_{21}' + E_{21} + G_{1} + E_{31})$ 
(of Dynkin type $D_{7} \oplus A_{12}$);

\par \vskip 5pt \noindent
$\Delta (5) = (E_{11}' + E_{12}' + F_{1} + E_{13}' + E_{13} + G_{3} + 
E_{23}) + (E_{33}' + E_{32}' + F_{3} + E_{31}' + E_{31} + 
G_{1} + E_{21} + E_{21}' + F_{2} + E_{22}' + E_{22} + G_{2})$ 
(of Dynkin type $D_{7} \oplus D_{12}$);

\par \vskip 5pt \noindent
$\Delta (6) = (E_{11}' + E_{12}' + E_{13}' + F_{1}) + (E_{33} + G_{3} + 
E_{23} + E_{23}' + F_{2} + E_{22}' + E_{22} + G_{2} + 
E_{32} + E_{32}' + F_{3} + E_{31}' + E_{31} + G_{1} + E_{21})$ 
(of Dynkin type $D_{4} \oplus A_{15}$).

\par \vskip 5pt
Let $S_{3} \rightarrow S(i)$ ($1 \leq i \leq 6$) be the contraction of 
$\Delta (i)$. Then the automorphism $g_{3}$ descends to automorphisms 
of $S(i)$, and the quotient surfaces 
$Z(i) := S(i)/<g_{3}>$ ($1 \leq i \leq 6$) 
are extremal log Enriques surfaces of type 
$D_{19}$, $D_{16} \oplus A_{3}$, 
$D_{13} \oplus A_{6}$, $D_{7} \oplus A_{12}$, $D_{7} \oplus D_{12}$ and 
$D_{4} \oplus A_{15}$ respectively.

\item 
Let $\Delta (7)$ be the divisor on $S_{2}$ defined by 

\par \vskip 5pt \noindent
$\Delta (7) = H_{31} + E_{31}' + F_{3} + E_{33}' + H_{33} + E_{33} + 
G_{3} + E_{13} + H_{13} + E_{13}' + F_{1} + E_{11}' + 
H_{11} + E_{11} + G_{1} + E_{21}' + F_{2} + E_{22}' + G_{2}$ 
(of Dynkin type $A_{19}$).

\par \vskip 5pt
Let $S_{2} \rightarrow S(7)$ be the contraction of $\Delta (7)$. 
Then the automorphism $g_{2}$ descends to an automorphism
of $S(7)$, and the quotient surface  $Z(7) := S(7)/<g_{2}>$  is an 
extremal log Enriques surface of type $A_{19}$. 
\endroster
\endproclaim 

\par
\demo{Proof} Since each connected component of $\text{Supp}\Delta (i)$ 
($1 \leq i \leq 6$) is $g_{3}-$stable, $g_{3}$ descends 
to its namesake on  $S(i)$. Since in addition every
1-dimensional component of  $S_{3}^{g_{3}}$
lies in  Supp $\Delta (i)$  and since no
connected component of  $\Delta(i)$  is disjoint from
$S_3^{g_3}$, it follows that the quotient map
$S(i) \rightarrow Z(i)$ has no ramification curves and that $Z(i)$ has no 
Du Val singular points. Thus, $Z(i)$ is a log Enriques surface whose 
canonical cover is equal to the quotient map  $S(i) \rightarrow Z(i)$. 
This implies the assertion (1). 
The verification of (2) is also similar. 
\qed 
\enddemo 

\head 
\S 2. Global canonical cover of an extremal log Enriques surface
\endhead 

\par \noindent
{\bf Note (2.1).} In this section, we let  $Z$  be an 
extremal log Enriques surface
of index  $I$, and we shall use the notation in the Introduction.

\par \vskip 1pc
The goal of this section is to show the following: 

\par
\proclaim{Proposition (2.2)}
\roster 
\item 
The canonical index $I$ is either $2$ or $3$. 
\item 
In the case where $I = 2$, $(S, g)$ is isomorphic to 
Shioda-Inose's pair $(S_{2}, g_{2})$ of discriminant $4$ and 
$Z$ is isomorphic to the extremal log Enriques surface $Z(7)$ 
defined in (1.7)(2). 
\item 
In the case where $I = 3$, $(S,g)$ is isomorphic to 
Shioda-Inose's pair $(S_{3}, g_{3})$ of discriminant $3$ and 
the type of $Z$ is either $D_{19}$, 
$D_{3l+1} \oplus D_{3m}$ or $D_{3l+1} \oplus A_{3m}$, where 
$l$ and $m$ are positive integers with $l + m = 6$.
\endroster
\endproclaim 
This Proposition will immediately follow from Lemmas 
(2.4), (2.8), (2.9), (2.11) and (2.13) below.  
First we remark some easy facts.

\par
\proclaim{Lemma (2.3)}
\roster
\item 
Every curve in $S^{[<g>]}$ is contained in $\Delta$. 
In particular, $S^{[<g>]}$ consists of smooth 
rational curves and finite isolated points.
\item 
$\Delta$ is $g-$stable, that is, $g(\Delta) = \Delta$.
\endroster

\par
\endproclaim
\demo{Proof} Since $K_{Z} \equiv 0$, the quotient map 
$\overline{S} \rightarrow Z$ is unramified in codimension one. 
This implies the assertion (1). The assertion (2) is clear.
\qed 
\enddemo
\proclaim{Lemma (2.4)} 
$I$ is either $2$, $3$, $4$, or $6$. 
\endproclaim
\demo{Proof} Since $S$ is a singular K3 surface, we know that 
$\text{rank}T_{S} = 2$, where $T_{S}$ denotes the transcendental 
lattice of $S$. Since $g^{*}\omega_{S} = \zeta_{I}\omega_{S}$ and 
$\omega_{S} \in T_{S}\otimes \Bbb C$, the action $g^{*}$ on $T_{S}$ 
has an eigen value $\zeta_{I}$. Thus, $\varphi(I) \leq 2$, where $\varphi$ is the Euler function. This implies the result. 
\qed 
\enddemo
We quote here the next two easy but useful Lemmas from [OZ1].

\par
\proclaim{Lemma (2.5) ([OZ1, Lemma 3.2])} 
Let $T$ be a smooth K3 surface and $\tau$ an involution of $T$ such that 
$\tau^{*}\omega_{T} = - \omega_{T}$.
\roster
\item 
Let $C_{1}$ and $C_{2}$ be two $g-$stable smooth rational curves 
on $T$ with $C_{1} \cdot C_{2} = 1$. Then, exactly one of $C_{i}$ is 
$\tau-$fixed.

\item 
Let $C$ be a $\tau-$stable but not $\tau-$fixed smooth rational 
curve on $T$. Set $C \cap T^{\tau} = \{P_{1}, P_{2}\}$.
Then, for each $i = 1, 2$, there exists a 
$\tau-$fixed curve $D_{i}$ passing through $P_{i}$. 
\endroster
\endproclaim

\par
\proclaim{Lemma (2.6) ([OZ1, Lemma 2.2, Proof of Lemma 2.3])} 
Let $T$ be a smooth K3 surface with an automorphism $\tau$ of $T$. 
Assume that $\tau$ is of order $3$ and that 
$\tau^{*}\omega_{T} = \zeta_{3}\omega_{T}$. 
\roster
\item 
Let $C_{1} + C_{2} + C_{3}$ be a linear chain of smooth rational curves 
on $T$. Assume that each $C_{i}$ is $\tau-$stable. Then, exactly 
one of $C_{i}$  is $\tau-$ fixed.

\item 
Let $C$ be a $\tau-$stable but not $\tau-$fixed smooth rational curve 
on $T$. Then, there exists a $\tau-$fixed curve $D$ on $T$ with 
$C \cdot D = 1$. 
\item 
Assume that $T^{\tau}$ consists of rational curves and isolated points. 
Let $N$ (resp. $M$) be the number of rational curves (resp. isolated points) 
in $T^{\tau}$. Then $M - N = 3$. 
\endroster
\endproclaim 
We return to our initial situation (2.1).

\par
\proclaim{Lemma (2.7)} 
Assume that $I = 2$.  Then we have:

\roster
\item 
Each connected component $\Delta_{i}$ of $\Delta$ is $g-$stable.
\item 
$\Delta_{i}$ is of type $A_{2n_{i}+1}$ ($n_{i} \in \Bbb Z_{\geq 0}$).
(See [Z1, Lemma 3.1].)

\item 
Each irreducible component of $\Delta_{i}$ is $g-$stable. 
\item   
Let $\Delta_{i} = C_{1} + C_{2} + \cdots + C_{2n_{i}+1}$ be 
the irreducible decomposition of $\Delta_{i}$ such that the dual graph 
of $\Delta_{i}$ is $C_{1} - C_{2} - \cdots - C_{2n_{i}+1}$. Then, 
$C_{j}$ is $g-$fixed if and only if  $j \equiv 1 (\text{mod} 2)$.

\endroster
\endproclaim

\par
\remark{Remark} This Lemma requires our assumption that a log Enriques 
surface has no Du Val singular points.
\endremark 

\par
\demo{Proof} We proceed the proof dividing into four steps.

\par
\proclaim{Step 1} Each $\Delta_{i}$ is $g-$stable. \endproclaim 

\par
\demo{Proof} This follows from our assumption that $Z$ has 
no Du Val singular points. 
\qed
\enddemo 

\par
\proclaim{Step 2} $\Delta_{i}$ is of type $A_{m}$ for certain integer $m$.
\endproclaim

\par
\demo{Proof} Assume the contrary that $\Delta_{i}$ is 
not of type $A_{m}$.
Then, by the classification of Dynkin diagram, 
there exists a unique rational curve $C$ 
in $\Delta_{i}$ which meets exactly three rational curves in 
$\Delta_{i}$, say, $D_{1}$, $D_{2}$, and $D_{3}$. 
Note that at least one of $D_{j}$, say $D_{1}$, does not 
meet any curves in 
$\Delta_{i}$ except for $C$.
By the uniqueness of $C$, we have $g(C) = C$ and 
$g(\{D_{1}, D_{2}, D_{3}\}) = \{D_{1}, D_{2}, D_{3}\}$. 
We shall derive a contradiction dividing into the two cases: \par
{\it Case 1.} $g \vert C = id$ and {\it Case 2.} $g \vert C \not= id$. \par
{\it Case 1.}  In this case, $D_{1}$ is $g-$stable but not $g-$fixed 
((2.5)(1)) and $D_1^{g}$ consists of two points. Since one of these two 
points is not in $C$, there exists a $g-$fixed curve $E (\not= C)$ 
which meets $D_{1}$ ((2.5)(2)).  This implies  $E \subset \Delta_{i}$ 
((2.3)(1)), a contradiction to the choice of $D_{1}$. 

\par 
{\it Case 2.} There exists exactly one $D_{j}$ with 
$g(D_{j}) = D_{j}$. Thus, $C^{g}$ contains a point $Q$ which does 
not lie in $D_{1} \cup D_{2} \cup D_{3}$. Then, there exists a $g-$fixed 
curve $E (\not= D_{1}, D_{2}, D_{3})$ passing through  $Q$ ((2.5)(2)). 
This implies  $E \subset \Delta_{i}$ ((2.3)(1)), a contradiction 
to the choice of $C$.
\qed 
\enddemo

\par
\proclaim{Step 3} 
Write  $\Delta_i = \sum_{j=1}^m C_j$
with  $C_j . C_{j+1} = 1$ ($1 \le j \le m-1$).
Then each irreducible component $C_{j}$  
is $g-$stable. 
\endproclaim

\demo{Proof} Assume the contrary that $g(C_{j}) \not= C_{j}$ 
for some $j$.  Then, $g(C_{j}) = C_{m + 1 - j}$ for all $j$, because 
$\text{\rm Aut}_{\text{\rm graph}}(A_{m}) \simeq \Bbb Z/2$. 
We shall drive a contradiction 
dividing into the two cases: \par 
{\it Case 1.} $m \equiv 0$ (mod 2), 
{\it Case 2.} $m \equiv 1$ (mod 2). 

\par 
{\it Case 1.} Since $g(C_{m/2} \cap C_{m/2+1}) = C_{m/2} \cap C_{m/2+1}$, 
there exists a $g-$fixed curve $E (\not= C_{m/2}, C_{m/2+1})$ 
passing through the point $C_{m/2} \cap C_{m/2+1}$. 
Then $E \subset \Delta_{i}$ ((2.3)(1)), a contradiction. 
\par 
{\it Case 2.} Since $C_{(m+1)/2}$ is $g-$stable but not $g-$fixed, 
there exists a $g-$fixed curve $E$ meeting $C_{(m+1)/2}$. Then 
$E \subset \Delta_{i}$, a contradiction.
\qed
\enddemo

\par
\proclaim{Step 4} 
$m \equiv 1 (\text{mod} 2)$; and $g \vert C_{j} = id$ if 
$j \equiv 1 (\text{mod} 2)$, 
while $g \vert C_{j}$ is an involution if $j \equiv 0 (\text{mod} 2)$.
\endproclaim 

\par
\demo{Proof} It follows from Step 3, (2.5)(1) and (2.3)(1) that 
both $C_{1}$ and $C_{m}$ are $g-$fixed. Now the result readily 
follows from (2.5)(1). 
\qed 
\enddemo
This completes the proof of (2.7). 
\qed 
\enddemo 

\par
\proclaim{Lemma (2.8)} 
Assume that $I = 2$. Then $(S, g)$ is isomorphic to 
Shioda-Inose's pair $(S_{2}, g_{2})$ of discriminant $4$ and 
$Z$ is isomorphic to the extremal log Enriques surface $Z(7)$ 
defined in (1.7)(2). 
\endproclaim
\demo{Proof} Let $N$ be the number of $g-$fixed curves on $S$. 
Recall that $\Delta_{i}$ contains just 
$(n_{i}+1)$ $g-$fixed curves (2.7) and that every $g-$fixed curve 
is contained in $\Delta = \coprod_{i=1}^{r} \Delta_{i}$ (2.3)(1). 
Thus, we get 
$$
N = \sum_{i = 1}^{r} (n_{i} + 1) = r + \sum_{i = 1}^{r} n_{i}. \tag 1
$$ 
On the other hand, since $Z$ is extremal, we have 
$$ 
19 = \text{rank}\Delta = \sum_{i = 1}^{r} (2n_{i}+1) 
= r + 2(\sum_{i = 1}^{r} n_{i}). \tag 2 
$$ 
Combining these two equalities, we get 
$$ 
N = r + (19 - r)/2 = (19 + r)/2 \geq 10. \tag 3 
$$ 
Now we may apply (1.6) to get $(S, g) \simeq (S_{2}, g_{2})$ and $N = 10$. 
This implies that $r = 1$ and that $\Delta$ is of type $A_{19}$. 
In other words, $Z$ is the most extremal log Enriques surface of 
type $A_{19}$. Now the result follows from  [OZ1, main Theorem 2]. 
\qed 
\enddemo 

\par
\proclaim{Lemma (2.9)} $I \not= 4$. \endproclaim 
\demo{Proof} Assume the contrary that $I = 4$. Then 
$h := g^{2}$ is an involution of $S$ with properties that 
$h^{*}\omega_{S} = -\omega_{S}$ and  $S^{h} \subset \Delta$. 
\proclaim{Claim} Each $\Delta_{i}$ is $h-$stable. \endproclaim 
\demo{Proof} Assume the contrary that 
$h(\Delta_{i}) \not= \Delta_{i}$ for some $i$. Then 
$\Delta_{i}$, $g(\Delta_{i})$, $h(\Delta_{i})$, and $g^{3}(\Delta_{i})$ 
are mutually different connected components of $\Delta$. Thus, 
$Z = \overline{S}/<g>$ has a Du Val singular point 
$\pi\circ\nu(\Delta_{i})$, a contradiction.
\qed
\enddemo 
By virtue of this Claim, we may repeat the same argument as in 
Steps 2-4 in (2.7) and (2.8) for the pair $(S, h)$  to show that 
$\Delta$ is of type $A_{19}$. This implies that $Z$ is the most extremal 
log Enriques surface of type $A_{19}$. However, then $I = 2$ 
by [OZ1, main Theorem 2], a contradiction. 
\qed
\enddemo 
\proclaim{Lemma (2.10)} Assume that $I = 3$. 
\roster
\item 
Each $\Delta_{i}$ is $g-$stable. Moreover, each irreducible 
component of $\Delta_{i}$ is also $g-$stable.
\item 
$\Delta_{i}$ is either of type $A_{n_{i}}$ or of type $D_{m_{i}}$ with 
$m_{i} \not\equiv 2 (\text{mod} 3)$.
(See [Z1, Prop. 6.1].)

\item 
All possible configurations of  $\Delta_{i}$
are given as follows, where we denote by $f$
(resp. $s$) $g-$fixed (resp. $g-$stable but not $g-$fixed)
irreducible components:

$$f \,\,\, \text{\rm (type} \, A_1), \,\,\,\, f-s \,\,\, 
\text{\rm (type} \, A_2), \,\,\,\,
s-f-s \,\,\, \text{\rm (type} \, A_3)$$

$$s-f-s-s-f-s-s- \dots -f-s-s-f-s \,\,\, \text{\rm (type} \, A_{3p})$$

$$f-s-s-f-s-s- \dots -f-s-s-f-s \,\,\, \text{\rm (type} \, A_{3p-1})$$

$$f-s-s-f-s-s- \dots -f-s-s-f \,\,\, \text{\rm (type} \, A_{3p-2})$$

\par \hskip 1.8pc
$s$
\par \hskip 1.8pc
$|$
\par
$$f-s-s-f-s-s- \dots -f-s-s-f-s \,\,\, \text{\rm (type} \, D_{3q+1})$$
\par \hskip 1.8pc
$|$
\par \hskip 1.8pc
$s$

\par \hskip 3pc
$s$
\par \hskip 3pc
$|$
\par
$$f-s-s-f-s-s- \dots -f-s-s-f \,\,\, \text{\rm (type} \, D_{3q})$$
\par \hskip 3pc
$|$
\par \hskip 3pc
$s$

\par
In particular, the pair of 
\newline 
(the number of $g-$fixed curves, the number of $g-$fixed isolated points) 
\newline
for each $\Delta_{i}$ is as follows: 

\par 
$(p, p+1)$ if of type $A_{3p}$; $(p, p)$ if of type $A_{3p-1}$; 
$(p, p-1)$ if of type $A_{3p-2}$; \par 
$(q, q+2)$ if of type $D_{3q+1}$; $(q, q+1)$ if of type $D_{3q}$.
\endroster
\endproclaim

\par
\demo{Proof} Making use of (2.3) and (2.6) (instead of (2.3) and (2.5))
and tracing out Dynkin diagrams,  
we can prove (2.10) in the same manner as in (2.7). 
Details will be left to the readers.
\qed
\enddemo

\par
\proclaim{Lemma (2.11)} 
Assume that $I = 3$. Then, $(S,g)$ is isomorphic to 
Shioda-Inose's pair $(S_{3}, g_{3})$ of discriminant $3$ and 
the type of $Z$ is either $D_{19}$, 
$D_{3l+1} \oplus D_{3m}$ or $D_{3l+1} \oplus A_{3m}$, where 
$l$ and $m$ are positive integers with $l + m = 6$.
\endproclaim 
\demo{Proof} Let $N$ (resp. $M$) be the number of $g-$fixed curves 
(resp. $g-$fixed isolated points) on $S$. 
Then by (2.6)(3), we have
$$ 
M - N = 3. \tag 1 
$$
On the other hand, we know by (2.10) that $\Delta$ 
is a disjoint sum of $a + b + c + d + e$ divisors whose types are: \par
$D_{3l_{1}+1}, \cdots , D_{3l_{a}+1}$, $D_{3m_{1}}, \cdots , D_{3m_{b}}$,
\par 
$A_{3p_{1}}, \cdots , A_{3p_{c}}$, 
$A_{3q_{1}-1}, \cdots , A_{3q_{d}-1}$, 
and $A_{3r_{1}-2}, \cdots , A_{3r_{e}-2}$, 

\par 
where $a$, $b$, $c$, $d$, and $e$ are certain non-negative integers. \par 
Then using (2.3)(1) and (2.10)(3), we calculate 
$$ 
N = \sum_{i = 1}^{a} l_{i} + \sum_{j = 1}^{b} m_{j} + 
\sum_{k = 1}^{c} p_{k} + \sum_{l = 1}^{d} q_{l} + \sum_{m = 1}^{e} r_{m} 
\tag 2
$$ 
and \par 
$
M \geq \sum_{i = 1}^{a} (l_{i}+2) + \sum_{j = 1}^{b} (m_{j}+1) + 
\sum_{k = 1}^{c} (p_{k}+1) + \sum_{l = 1}^{d} q_{l} + 
\sum_{m = 1}^{e} (r_{m} - 1) 
$
\par 
$$
= \sum_{i = 1}^{a} l_{i} + \sum_{j = 1}^{b} m_{j} + 
\sum_{k = 1}^{c} p_{k} +
\sum_{l = 1}^{d} q_{l} + \sum_{m = 1}^{e} r_{m} 
+ (2a + b + c -e). \tag 3 
$$
\par 
Substituting (2) and (3) into (1), we get 
$$ 
3 = M - N \geq 2a + b + c - e. \tag 4 
$$ 
Since $Z$ is an extremal log Enriques surface, we calculate 
\par 
$19 = \sum_{i = 1}^{a} (3l_{i}+1) + \sum_{j = 1}^{b} 3m_{j} + 
\sum_{k = 1}^{c} 3p_{k} + \sum_{l = 1}^{d} (3q_{l}-1) + 
\sum_{m = 1}^{e} (3r_{m} - 2)$
\par 
$= 3N + a - d - 2e$, 
\par 
where we use (2) to get the last equality. Thus, 
$$ 
N = \frac{19 + (2e + d - a)}{3}. \tag 5 
$$ 
Suppose that $2e + d - a \leq -2$. Then $a \geq 2e + d +2$. 
Substituting this into (4), we get 
$$ 
3 \geq 2a + b + c - e \geq 2(2e + d + 2) + b + c - e 
= 4 + b + c + 2d + 3e \geq 4,  
$$ 
a contradiction. Thus $2e + d - a \geq -1$. Substituting this into (5), 
we get $N \geq (19 - 1)/3 = 6$. Now we may apply (1.3) to get 
$(S, g) \simeq (S_{3}, g_{3})$ and then $N = 6$. 
Combining this equality with (5), we get $2e + d - a = -1$, that is, 
$a = 2e + d + 1$. Substituting this into (4), we calculate
$$ 
3 \geq 2(2e + d + 1) + b + c - e = 2 + 2d + b + c + 3e, 
$$ 
that is, 
$$ 
1 \geq 2d + 3e + b + c. 
$$ 
>From this, we can easily see that $d = e = 0$, 
$a = 2e + d + 1 = 1$ and $b + c \leq 1$. 
Combining these formula with $\text{rank}\Delta =19$, we readily 
see that $\Delta$ is either one of the following types: 
$D_{19}$, $D_{3l+1} \oplus D_{3m}$, or $D_{3l+1} \oplus A_{3m}$ 
($l + m = 6$). This is nothing but the last half assertion of (2.11). 
\qed 
\enddemo 

\par
It only remains to show $I \not= 6$. For this we need the following:

\par
\proclaim{Proposition (2.12)} 
Let $T$ be a smooth K3 surface and $\tau$ an automorphism of $T$. 
Assume that 
\roster
\item
$\tau$ is of order $6$ and $\tau^{*}\omega_{T} = \zeta_{6}\omega_{T}$, 
and that
\item 
$T^{[<\tau>]}$ consists only of isolated points and smooth rational curves.
\endroster
Then $T^{\tau} (= T^{\tau^{5}})$, $T^{\tau^{2}} (= T^{\tau^{4}})$, and 
$T^{\tau^{3}}$ are described as follows: 
$$ 
T^{\tau} = (\coprod_{i=1}^{2(c+1)}\{P_{i}\}) \coprod (\coprod_{i=1}^{2(c+1)} 
\{Q_{i}\}) \coprod (\coprod_{j=1}^{c} 
C_{j}), 
$$ 
$$ 
T^{\tau^{2}} = (\coprod_{i=1}^{2(c+1)}\{P_{i}\}) 
\coprod (\coprod_{k=1}^{2(p+1)} \{P_{k}'\}) 
\coprod (\coprod_{j=1}^{c} C_{j}) 
\coprod (\coprod_{l=1}^{c+1} D_{l}) 
\coprod (\coprod_{m=1}^{2p} F_{m}),
$$
$$
T^{\tau^{3}} = (\coprod_{j=1}^{c} C_{j})  
\coprod (\coprod_{i=1}^{2(c+1)} E_{i}) 
\coprod (\coprod_{n=1}^{3q} G_{n}),
$$ 
where $c$, $p$, and $q$ are non-negative integers with 
$c + p + q \leq 2$, $P_{*}$, $Q_{*}$, and $P_{*}'$ are isolated
points, and  $C_{*}$, $D_{*}$, $E_{*}$, $F_{*}$, $G_{*}$ 
are smooth rational curves. Moreover, each of $D_{*}$ 
and $E_{*}$ is $\tau-$stable, while 
$\tau$ acts on each set $\{F_{2i-1}, F_{2i} \}$ as an involution and  
on $\{G_{3j-2}, G_{3j-1}, G_{3j} \}$ of order $3$.
\endproclaim 

\par
\demo{Proof} Our proof is based on the holomorphic 
Lefschetz fixed point formula [AS1, 2], the topological Lefschetz fixed 
point formula [U], and local coordinate calculation. 
We shall divide the proof into three steps. 
\proclaim{Step 1} \par 
$T^{\tau} = \{P_{1}\} \coprod \cdots \coprod \{P_{2l}\} \coprod 
\{Q_{1}\} \coprod \cdots \coprod \{Q_{2l}\} \coprod 
C_{1} \coprod \cdots \coprod C_{c}$, \par 
$T^{\tau^{2}} = \{P_{1}\} \coprod \cdots \coprod \{P_{2l}\} 
\coprod \{P_{1}'\} \coprod \cdots \coprod \{P_{k'}'\}$ \par 
$\coprod C_{1} \coprod \cdots \coprod C_{c} 
\coprod D_{1} \coprod \cdots \coprod D_{l} 
\coprod F_{1} \coprod \cdots \coprod F_{p'}$, \par
$T^{\tau^{3}} = C_{1} \coprod \cdots \coprod C_{c} 
\coprod E_{1} \coprod \cdots \coprod E_{2l} 
\coprod G_{1} \coprod \cdots \coprod G_{q'}$, \par 
where $l$, $c$, $p'$, $k$ and $q'$ are non-negative integers and 
$C_{*}$, $D_{*}$, $F_{*}$, $E_{*}$, and $G_{*}$ are 
smooth rational curves. 
Moreover $Q_{2i-1}, Q_{2i} \in D_{i}$, $P_{j}, Q_{j} \in E_{j}$, and 
each of $D_{i}$ and $E_{j}$ is $\tau-$stable.
\endproclaim 

\par
\demo{Proof} Suppose that $P$ is an isolated point of $T^{\tau}$. 
Since $\tau^{*}\omega_{T} = \zeta_{6}\omega_{T}$, there exist 
local coordinates $(x_{P}, y_{P})$ around $P$ such that either
\roster
\item 
$\tau^{*}(x_{P}, y_{P}) = (\zeta_{6}^{2}x_{P}, \zeta_{6}^{5}y_{P})$ 
or 
\item 
$\tau^{*}(x_{P}, y_{P}) = (\zeta_{6}^{3}x_{P}, \zeta_{6}^{4}y_{P})$.
\endroster 
Denote by $P_{1}, \cdots , P_{a} (\in T^{\tau})$ (resp. 
by $Q_{1}, \cdots , Q_{b} (\in T^{\tau}))$ the points of type (1) (resp. 
of type (2)). Then we write 
$T^{\tau} =  \{P_{1}\} \coprod \cdots \coprod \{P_{a}\} \coprod 
\{Q_{1}\} \coprod \cdots \coprod \{Q_{b}\} \coprod
C_{1} \coprod \cdots \coprod C_{c}$, where $C_{\alpha}$ 
are smooth rational curves. Let $R$ be a point in $C_{\alpha}$. 
Then there exist local 
coordinates $(x_{R}, y_{R})$ around $R$ such that 
$\tau^{*}(x_{R}, y_{R}) = (x_{R}, \zeta_{6}y_{R})$.
Note that  $C_{\alpha} = (y_{R} = 0)$ around $R$. 

\par 
Let $P$ be a point in $\{P_{1}, ... ,P_{a}\}$. 
Since $(\tau^{*})^{2}(x_{P}, y_{P}) = 
(\zeta_{6}^{4}x_{P}, \zeta_{6}^{4}y_{P})$ by (1), $P$ is an 
isolated $\tau^{2}-$fixed point. 

\par 
Let $Q$ be a point in $\{Q_{1}, ... ,Q_{b}\}$. 
Since $(\tau^{*})^{2}(x_{Q}, y_{Q}) = 
(x_{Q}, -y_{Q})$ by (2), there exists a unique smooth rational curve, 
say $D$, such that $Q \in D \subset T^{\tau^{2}}$.  
Note that  $D = (y_{Q} = 0)$ around $Q$. In particular, 
$D$ is $\tau-$stable and 
$\tau \vert D$ is an involution on $D$. Thus $\tau$ has another 
$\tau-$fixed point $Q'$ on $D$ 
around which $(\tau\vert D)^{*} = (\zeta_{6}^{4})^{-1}$. 
Since $C_{\alpha}$ and $D$ are 
disjoint (by the smoothness of $T^{\tau^{2}}$), this point $Q'$ 
is also isolated in  $T^{\tau}$  and in fact contained
in $\{Q_{1}, ... ,Q_{b}\}$.  
Now setting $\{D_{j} (1 \leq j \leq l) \vert D_{j} \subset T^{\tau^{2}}, 
D_{j} \cap \{Q_{1}, ... ,Q_{b}\} \not= \phi \}$, and using the smoothness 
of $T^{\tau^{2}}$, we can adjust the numbering of $Q_{\beta}$ 
($1 \leq \beta \leq b$) as 
$Q_{1}, Q_{2} \in D_{1}, ... , Q_{b-1}, Q_{b} \in D_{l}$. 
In particular, $b = 2l$. 

\par 
Next we examine $T^{\tau^{3}}$. Again, let $P$ (resp. $Q$) be a 
point in $\{P_{1}, ... ,P_{a}\}$ (resp. in $\{Q_{1}, ... ,Q_{2l}\}$). 
Since $(\tau^{*})^{3}(x_{P}, y_{P}) = 
(x_{P}, -y_{P})$ around $P$, there exists a unique smooth rational 
curve $E'$ such that $P \in E' \subset T^{\tau^{3}}$ (and that 
$E' = (y_{P} = 0)$ around $P$). Similarly, there exists a unique smooth
rational curve $E''$ such that $Q \in E'' \subset T^{\tau^{3}}$ (and that 
$E'' = (x_{Q} = 0)$ around $Q$). 
Using this description, we easily see that 
both $E'$ and $E''$ are $\tau-$stable and that 
$\tau \vert E'$ is a multiplication 
by $\zeta_{6}^{2}$ around $P$ and $\tau \vert E''$ is a multiplication 
by $\zeta_{6}^{4}$ around $Q$. Note also that $\vert(E')^{\tau}\vert = 2$ and 
$\vert(E'')^{\tau}\vert = 2$. 

\par 
Let $E_{i}$ ($1 \leq i \leq m$) be the $\tau^{3}-$fixed curves which 
contains at least one point in 
$\{P_{\alpha}, Q_{\beta} \vert 1 \leq \alpha \leq a, 1 \leq \beta \leq 2l\}$. By the smoothness of $T^{\tau^{3}}$, each $E_{i}$ coincides with 
some $E'$ or $E''$ found in the above process. In particular, each $E_{i}$ is $\tau-$stable. Then, using again the smoothness of $T^{\tau^{3}}$
and the description of $T^{\tau}$, and regarding the two points 
$E^{\tau}$ as $0$ and $\infty$ of $E_{i} (\simeq \Bbb P^{1})$, we see that 
there exist bijections $\varphi : \{1, ... , m\} \rightarrow \{1, ... , a\}$ 
and $\psi : \{1, ... , m\} \rightarrow \{1, ... , 2l\}$ such that 
$E_{i}^{\tau} = \{P_{\varphi(i)}, Q_{\psi(i)}\}$. Thus $m = a = 2l$. Then renumbering 
$E_{*}$ and $P_{*}$, we have $P_{i}, Q_{i} \in E_{i}$ for all $i$ with 
$1 \leq i \leq 2l$. Since $T^{\tau^{3}}$ contains no isolated points, we can now easily get the description of $T^{\tau^{3}}$ in Step 1.
Now we get the desired description of  $T^{\tau}, T^{\tau^2},
T^{\tau^3}$.
\qed 
\enddemo

\par
\proclaim{Step 2} 
$l = c + 1$, $p' = 2p$, $q' = 3q$ and $k = 2(p+1)$ for some non-negative 
integers $p$ and $q$, where $l$, $p'$, $q'$ and $k$ are integers found 
in Step 1.
\endproclaim

\par
\demo{Proof} We apply the holomorphic Lefschetz fixed point formula [AS1,2]
for  $(T, \tau)$:
$$ 
L(\tau) := \sum (-1)^{i}\text{tr}(\tau^{*} \vert H^{i}(T, \Cal O_{T})) 
= \sum_{j=1}^{2l} a(P_{j}) + \sum_{j=1}^{2l} a(Q_{j}) 
+ \sum_{i=1}^{c} b(C_{i}).
$$ 
We calculate both sides and compare them. 

\par 
Using the Serre duality, we get from the first equality that 
$$ 
L(\tau) = 1 + \zeta_{6}^{-1} = \frac{3 - \sqrt{-3}}{2}. \tag 1
$$
By the definition of $a(*)$ aa  in [AS1,2] 
and the local description of  $\tau$-action given in Step 1, 
we calculate 

\par
$$a(P_{j}) := \frac{1}{det(1 - \tau^{*}\vert T_{P_{j}})} 
= \frac{1}{(1 - \zeta_{6}^{2})(1 - \zeta_{6}^{5})} 
= \frac{3 - \sqrt{-3}}{6},$$ 
$$a(Q_{j}) = \frac{1}{(1 - \zeta_{6}^{3})(1 - \zeta_{6}^{4})} 
= \frac{3 - \sqrt{-3}}{12},$$
$$b(C_{i}) := \frac{1 - g(C_{i})}{1 - \zeta_{6}} 
- \frac{\zeta_{6}}{(1 - \zeta_{6})^{2}} \cdot (C_{i}^{2}) 
= \frac{1}{1 - \zeta_{6}} 
- \frac{\zeta_{6}}{(1 - \zeta_{6})^{2}} \cdot (-2) 
= \frac{-(3 - \sqrt{-3})}{2}.$$ 

\par
Using the above formula for  $L(\sigma)$  in terms 
of  $a(*), b(*)$, we obtain:
$$
L(\tau) = \frac{3 - \sqrt{-3}}{6} \cdot 2l 
+ \frac{3 - \sqrt{-3}}{12} \cdot 2l - \frac{(3 - \sqrt{-3})}{2} \cdot c. 
\tag 2 
$$ 
Combining (1) and (2), we readily see that $l = c + 1$. Thus, \par 
$T^{\tau^{2}} = \{P_{1}\} \coprod \cdots \coprod \{P_{2(c+1)}\} 
\coprod \{P_{1}'\} \coprod \cdots \coprod \{P_{k'}'\}$ \par 
$\coprod C_{1} \coprod \cdots \coprod C_{c} 
\coprod D_{1} \coprod \cdots \coprod D_{(c+1)} 
\coprod F_{1} \coprod \cdots \coprod F_{p'}$. \par
Using this description and the smoothness of $T^{\tau^{2}}$, 
we easily see that $\tau$ acts on both $\{F_{1}, ... , F_{p'}\}$ 
and $\{P_{1}', ... , P_{k'}'\}$ as fixed point free involutions. 
Thus, $p' = 2p$ and $k' = 2k$ for some integers $p$ and $k$. 

\par
Next, we shall find a relation between  $k$  and  $p$. 
Applying (2.6) to the pair  $(T, \tau^2)$  where
ord$(\tau^2) = 3$, we obtain
$\#(\tau^{2}-\text{isolated points}) - \#(\tau^{2}-\text{fixed curves}) = 3$, 
that is, 

\par 
$2(c+1) + 2k - (c + (c + 1) + 2p ) = 3$. 

\par 
This implies $k = p + 1$. 
Using the description of $T^{\tau^{3}}$ and 
applying the same argument as before for the set
$\{G_1, G_2, \dots, G_{q'}\}$ (instead of  $\{F_1, F_2, \dots, F_{p'}\}$),
we can readily see that 
$\tau$ induces a fixed point free automorphism of order 3 
on the set $\{G_{1}, ... , G_{q'}\}$. Thus, $q' = 3q$ for some integer $q$. 
This completes Step 2. 
\qed 
\enddemo 

\par
Now we only remain to show the inequality $c + p + q \leq 2$. 
Let us consider the action $\tau^{*}$ on $H^{2}(T, \Bbb Q)$. 
Since $(\tau^{*})^{6} = id.$ and $\tau^{*}\omega_{T} = \zeta_{6}\omega_{T}$, 
the pairs of 
\newline
(the eigenvalue of  $\, \tau^{*} \vert H^{2}(T, \Bbb Q),$ 
its multiplicity) 
\newline
are written as
$$ 
(1, \alpha), (-1, \beta), 
(\zeta_{3}, \gamma), (\overline{\zeta_{3}}, \gamma), 
(\zeta_{6}, 1 + \delta), (\overline{\zeta_{6}}, 1 + \delta), 
$$ 
where $\alpha, \beta, \gamma$ and $\delta$ are certain non-negative integers. 
Now the required inequality follows from $\delta \geq 0$ and the next Step 3.

\par
\proclaim{Step 3} With the above notation, \par 
$\alpha = 5c + 2p + q + 6$, \par
$\beta = -c + 2p - q + 4$, \par 
$\gamma = -c - p + q + 3$, and \par 
$\delta = -c - p -q + 2$. 
\endproclaim 

\par
\demo{Proof} Since $\text{dim}H^{2}(T, \Bbb Q) = 22$, we have 
$$ 
\alpha + \beta + 2\gamma + 2\delta = 20. \tag 1
$$
In order to obtain other relations, we make use of the topological 
Lefschetz fixed point formula [U]:
$$
\chi_{top}(T^{\tau^j}) =  
\sum_{i = 0}^{4} (-1)^{i} \text{tr} ((\tau^{*})^{j} \vert H^{i}(S, \Bbb Q)). \tag *
$$ 
Using
$T^{\tau} = \{P_{1}\} \coprod \cdots \coprod \{P_{2(c+1)}\} \coprod 
\{Q_{1}\} \coprod \cdots \coprod \{Q_{2(c+1)}\} \coprod 
C_{1} \coprod \cdots \coprod C_{c}$ and applying (*) with  $j = 1$,
we get $4(c +1) + 2c = 2 + \alpha - \beta - \gamma + \delta + 1$. 
This gives 
$$
\alpha - \beta - \gamma + \delta = 6c + 1. \tag 2
$$
Next using 
$T^{\tau^{2}} = \{P_{1}\} \coprod \cdots \coprod \{P_{2(c+1)}\} 
\coprod \{P_{1}'\} \coprod \cdots \coprod \{P_{2(p+1)}'\}
\coprod C_{1} \coprod \cdots \coprod C_{c}$
\par
$\coprod D_{1} \coprod \cdots \coprod D_{c+1} 
\coprod F_{1} \coprod \cdots \coprod F_{2p}$  and applying (*) with  $j = 2$, we get 
\par
$2(c +1) + 2(p + 1) + 2c + 2(c+1) + 2\cdot 2p = 2 + (\alpha + \beta) - (\gamma + \delta + 1).$ 
\par 
This gives 
$$
\alpha + \beta - \gamma - \delta = 6c + 6p + 5. \tag 3
$$
Finally using
$T^{\tau^{3}} = C_{1} \coprod \cdots \coprod C_{c} 
\coprod E_{1} \coprod \cdots \coprod E_{2(c+1)} 
\coprod G_{1} \coprod \cdots \coprod G_{3q}$ and applying (*) for 
$\tau^{3}$, we get 
$2c + 2 \cdot 2(c+1) + 2 \cdot 3q = 2 + (\alpha + 2\gamma) - 
(\beta + 2(\delta + 1))$. This implies 
$$ 
\alpha - \beta + 2\gamma - 2\delta = 6c + 6q + 4. \tag 4
$$
Now solving the equations (1) - (4) for $\alpha$, $\beta$, 
$\gamma$, $\delta$, we get the result.
\qed
\enddemo
This completes the proof of (2.12). 
\qed 
\enddemo 

\par
Returning back to our intial setting (2.1), we prove the following:

\par
\proclaim{Lemma (2.13)} $I \not= 6$. 
\endproclaim 

\par
\demo{Proof} Assume that $I = 6$. Then applying (2.12) for 
$(S, g)$, we see that there are non-negative integers
$c, p, q$  such that  $c+p+q \le 2$  and
that the irreducible decompositions of the 1-dimensional 
locus of $S^{g}$, $S^{g^{2}}$ and $S^{g^{3}}$ are written as follows 
respectively: 

\par 
$C_{1} \coprod \cdots \coprod C_{c}$;

\par 
$C_{1} \coprod \cdots \coprod C_{c} 
\coprod D_{1} \coprod \cdots \coprod D_{c+1} 
\coprod F_{1} \coprod \cdots \coprod F_{2p}$; 

\par 
$C_{1} \coprod \cdots \coprod C_{c} 
\coprod E_{1} \coprod \cdots \coprod E_{2(c+1)} 
\coprod G_{1} \coprod \cdots \coprod G_{3q}$,

\par 
where $C_{*}$, $D_{*}$, and $E_{*}$ are $g-$stable 
while the other $F_{*}$ and $G_{*}$ are not $g-$stable. 
Note also by (2.3)(1) that these curves 
$C_*, D_*, E_*$  are all contained in $\Delta$. 

\par 
Let us consider the connected components $\Delta_{i}$  of $\Delta$. 
Since $Z = \overline{S}/<g>$ has no Du Val singular points, each $\Delta_{i}$ 
satisfies either
\roster
\item 
$g^{3}-$stable or 
\item 
$g^{2}-$stable but not $g-$stable.
\endroster
Let $\Delta_{i}$ ($1 \leq i \leq n$) be of type (1) and 
$\Delta_{j}$ ($n + 1 \leq j \leq n + m = r$) of type (2). 
Since $g^{3}$ is of order 2 and $(g^{3})^{*}\omega_{S} = -\omega_{S}$, 
it follows from the argument in (2.7) (Steps 2-4) that 
each $\Delta_{i}$ ($1 \leq i \leq n$) is of the Dynkin type $A_{2\alpha_{i}+1}$ and contains exactly $(\alpha_{i} + 1)$ $g^{3}-$fixed 
curves. On the other hand, the above description of $S^{g^{3}}$ 
shows that the number of all the $g^{3}-$fixed curves is just 
$3(c + q) + 2$. Thus, 
$$ 
\sum_{i = 1}^{n} \text{rank}\Delta_{i} 
= \sum_{i = 1}^{n} (2(\alpha_{i} + 1) - 1) 
= 2 \sum_{i = 1}^{n} (\alpha_{i} + 1) - n = 6(c + q) + 4 - n. \tag 1 
$$ 
Let us consider the connected components $\Delta_{j}$ of type (2). 
Since $g^{2}$ is of order 3 and $(g^{2})^{*}\omega_{S} = \zeta_{3}\omega_{S}$, 
it follows from the argument in (2.10) that 
each $\Delta_{j}$ is of Dynkin type $A_{*}$ 
or $D_{*}$ and contains at least one $g^{2}-$fixed curve. Moreover, only 
$F_{*}$ are the $g^{2}-$fixed curves in $\Delta_{j}$, because $C_{*}$ 
and $D_{*}$ are $g-$stable so they are in $\Delta_{i}$ ($1 \leq i \leq n$). 
Thus, 
$$ 
m \leq 2p, \tag 2 
$$
and  $n \ge 1$ (because there is at least one  $D_*$) 
and  $\text{rank}\Delta_{j} \leq 
3 \cdot \vert \{F_{*} | F_{*} \subset \Delta_{j}\} \vert + 1$ (2.10(3)). 
Thus, 
$$
\sum_{j = n +1}^{n + m} \text{rank}\Delta_{j} 
\leq 3\cdot 2p + m = 6p + m. \tag 3
$$ 
Combining (1), (2) and (3) with
$19 = \sum_{i = 1}^{n} \text{rank}\Delta_{i} 
+ \sum_{j = n +1}^{n + m} \text{rank}\Delta_{j}$ and $c + p + q \leq 2$, 
we get 
$$ 
19 \leq 6(c + q) + 4 - n + 6p + m = 6(c + p + q) + 4 + m - n 
\leq 6 \cdot 2 + 4 + 2p - n \leq 19. \tag 4
$$
Thus the all inequalities in (4) must be equalities. 
This implies $n = 1$, $p= 2$, $m = 2p = 4$, and $c = q = 0$. 
Combining these equalities with rank  $\Delta = 19$ , we readily see that 
$\Delta$ is of type $A_{3} \oplus D_{4}^{\oplus 4}$. 
Then using (2.10), we see that $\Delta$ contains $2 + 4\cdot3 = 14$ isolated 
$g^{2}-$fixed points and that $g^{2}$ has exactly 5 fixed curves. 
Thus, $M \geq 14$ and $N = 5$, where $M$ is the number of the 
isolated $g^{2}-$fixed points and $N$ is that of the $g^{2}-$fixed curves
on  $S$.  However this contradicts the equality 
$M - N = 3$ ((2.6)(3)). This completes the proof. 
\qed
\enddemo

\par
\head 
\S 3. Types of extremal log Enriques surfaces
\endhead

\par
The goal of this section is to finish the proof of the Main Theorem (1). 

\par
Let  $Z$  be an extremal log Enriques surface of index  $I$
and we shall use the notation in the Introduction.
By (2.11), we already know that $\Delta$ is either one of the following 
types: $A_{19}$, $D_{19}$, $D_{16} \oplus A_{3}$, $D_{13} \oplus A_{6}$, 
$D_{10} \oplus A_{9}$, $D_{7} \oplus A_{12}$, 
$D_{4} \oplus A_{15}$, $D_{13} \oplus D_{6}$, $D_{10} \oplus D_{9}$,
$D_{7} \oplus D_{12}$, $D_{4} \oplus D_{15}$. \par 
Thus, in order to get the Main Theorem (1), we may prove the following:

\par
\proclaim{Lemma (3.1)} 
\roster
\item 
$\Delta$ is not of types $D_{13} \oplus D_{6}$, $D_{10} \oplus D_{9}$, 
$D_{4} \oplus D_{15}$.
\item 
$\Delta$ is not of type $D_{10} \oplus A_{9}$.
\endroster
\endproclaim 

\par
\demo{Proof of (1)} We shall argue by contradiction. 
Since $(S, g) \simeq (S_{3}, g_{3})$ by (2.2), we may identify $(S, g)$ 
with $(S_{3}, g_{3})$. We denote 
$\text{Supp} \Delta = (\cup_{i=1}^{3l+1} C_{i}) \coprod 
(\cup_{j=1}^{3m} E_{j})$ where the numberings are given as 
$C_{1}.C_{3} = C_{2}.C_{3} = C_{i}.C_{i+1} = 1 (i \geq 3)$ and 
$E_{1}.E_{3} = E_{2}.E_{3} = E_{j}.E_{j+1} = 1 (j \geq 3)$.
By $\Delta$, we also denote the sublattice of 
$\text{Pic} S_{3}$ generated by the irreducible components of $\Delta$. 
Let us consider the primitive closure $\overline{\Delta}$ of 
$\Delta$ in $\text{Pic} S_{3}$. 
Since $[\overline{\Delta} : \Delta]^{2} = 
({\text{discr} \Delta})/({\text{discr} \overline{\Delta}}) = 
16/({\text{discr} \overline{\Delta}})$, $[\overline{\Delta} : \Delta]$ 
is either 1, 2 or 4. Dividing into these three cases, 
we shall derive a contradiction. 

\par
First assume that $[\overline{\Delta} : \Delta] = 4$. 
Then $\text{discr} \overline{\Delta} = 1$. Thus, we have an othorgonal decomposition 
of $\text{Pic} S_{3}$: $\text{Pic} S_{3} = \overline{\Delta} \oplus \Bbb Z \cdot H$. 
This implies $H^{2} = \text{discr Pic} S_{3} = 3 \not\equiv 0 (\text{mod} 2)$, 
a contradiction. \par
Next assume that $[\overline{\Delta} : \Delta] = 2$. 
Then there exist integers $\alpha_{i}, \beta_{j} \in \{0,1\}$ such that 
$L:= \frac{1}{2}(\sum_{i=1}^{3l+1} \alpha_{i}C_{i} 
+ \sum_{j = 1}^{3m} \beta_{j}E_{j}) \in \overline{\Delta} - \Delta$. 
Substituting $i = 3l+1, 3l, ..., 3$ and $j = 3m, 3m-1, ..., 3$ into 
$L.C_{i} \in \Bbb Z$ and $L.E_{j} \in \Bbb Z$, we readily find that \par 
$L =$ $\frac{1}{2}(E_{6} + E_{4} + E_{2})$ or $\frac{1}{2}(E_{6} + E_{4} + E_{2})$ in the case where $\Delta$ is of type $D_{13} \oplus D_{6}$, \par 
$L =$ $\frac{1}{2}(C_{10} + C_{8} + C_{6} + C_{4} + C_{2})$ or $\frac{1}{2}(C_{10} + C_{8} + C_{6} + C_{4} + C_{1})$ 
in the case where $\Delta$ is of type $D_{10} \oplus D_{9}$, and 

\par 
$L =$ $\frac{1}{2}(E_{4} + E_{2})$ or $\frac{1}{2}(E_{4} + E_{1})$ in the case where $\Delta$ is of type $D_{4} \oplus D_{15}$. \par 
But this contradicts the next Lemma due to Nikulin [N]:

\par
\proclaim{Lemma (3.2)} 
Let $C_{1}, C_{2}, ..., C_{l}$ be mutually disjoint smooth rational 
curves on a smooth K3 surface $T$. 
Assume that $C_{1} + \cdots + C_{l} \in 2 \cdot \text{Pic} T$. 
Then, $l$ is either $0$, $8$, or $16$. 
\qed \endproclaim

\par
Finally assume that $[\overline{\Delta} : \Delta] = 1$. That is, $\Delta$ 
is primitive in $\text{Pic} S_{3}$. Then there exists an element $h \in 
\text{Pic} S_{3}$ 
such that 
$\text{Pic} S_{3} = \Bbb Z <C_{1}, ..., C_{3l+1}, E_{1}, ..., E_{3m}, h>$. 
Set $\Bbb Z \cdot H = \Delta^{\bot}$ in $\text{Pic} S_{3}$ and 
$n = [\text{Pic} S_{3} : \Delta \oplus \Bbb Z \cdot H]$. 
Then $n^{2} = (\text{discr}\Delta \cdot H^{2})/(\text{discr Pic} S_{3}) = 
16 H^{2}/3$, that is, $H^{2} = \frac{3}{16}n^{2}$. 
On the other hand, by replacing $h$ by $-h$ if necessary, we can find 
integers $a_{i}$, $b_{j}$ such that 
$H = nh + \sum a_{i}C_{i} + \sum b_{j}E_{j}$ in $\text{Pic} S_{3}$, that is, 
$\frac{H}{n} = h + \sum \frac{a_{i}}{n}C_{i} + \sum \frac{b_{j}}{n}E_{j}$ 
in $\text{Pic} S_{3} \otimes \Bbb Q$. Using $H.C_{\alpha} = H.E_{\beta} = 0$ and 
the negative definiteness of $(C_{i}\cdot C_{\alpha})$ and 
$(E_{j}\cdot E_{\beta})$, we see that $(\frac{a_{i}}{n}, \frac{b_{j}}{n})$ 
is the unique solution of 
$$
(h + \sum x_{i}C_{i} + \sum y_{j}E_{j})\cdot C_{\alpha} = 0, 
(h + \sum x_{i}C_{i} + \sum y_{j}E_{j})\cdot E_{\beta} = 0, 
$$ 
($\alpha = 1, ..., 3l+1; \beta = 1, ..., 3m$). \par 
Since $\text{discr}(C_{i}\cdot C_{\alpha}) = \text{discr}(E_{j}\cdot E_{\beta}) = 4$, 
this implies that $\frac{a_{i}}{n}, \frac{b_{j}}{n} \in \frac{\Bbb Z}{4}$. 
Thus, 
$\frac{4H}{n} = 4h + \sum \frac{4a_{i}}{n}C_{i} 
+ \sum \frac{4b_{j}}{n}E_{j} \in \text{Pic} S_{3}$. 
This implies $n \vert4$. However, then $H^{2} = \frac{3}{16}n^{2} \not\equiv 0 (\text{mod} 2)$, a contradiction. This proves the assertion (1). 
\qed
\enddemo

\par
\demo{Proof of (2)} The verification of (2) is quite similar to 
that of (1). Assuming the contrary, we identify $(S, g)$ with 
$(S_{3}, g_{3})$ and set 
$\text{Supp} \Delta = (\cup_{i=1}^{10} C_{i}) \coprod 
(\cup_{j=1}^{9} E_{j})$ where the names are given as:
$\sum C_{i}$ is of type $D_{10}$, $C_{1}.C_{3} = C_{2}.C_{3} = C_{i}.C_{i+1} = 1 (i \geq 3)$, $\sum E_{j}$ is of type $A_{9}$ and  
$E_{j}.E_{j+1} = 1 (j \geq 1)$.
Let $\overline{\Delta}$ be the primitive closure of the sublattice 
$\Delta$ in $\text{Pic} S_{3}$. 
Again, it follows from $[\overline{\Delta} : \Delta]^2 = 
(\text{discr} \Delta)/(\text{discr} \overline{\Delta}) = 
40/(\text{discr} \overline{\Delta})$ that $[\overline{\Delta} : \Delta]$ 
is either 1 or 2. In each case, we shall derive a contradiction. 

\par 
First assume that $[\overline{\Delta} : \Delta] = 2$. 
Then there exists integers $\alpha_{i}, \beta_{j} \in \{0, 1\}$ 
such that 
$L := \frac{1}{2}(\sum \alpha_{i}C_{i} + \sum \beta_{j}E_{j}) \in 
\overline{\Delta} - \Delta$. Using $L\cdot C_{i}, L\cdot E_{j} \in \Bbb Z$, 
we find that $L$ is either one of  
$\frac{\delta_{1}}{2}(C_{10} + C_{8} + C_{6} + C_{4} + C_{i}) + 
\frac{\delta_{2}}{2}(E_{9} + E_{7} + E_{5} + E_{3} + E_{1})$, 
where $i \in \{1, 2\}$, $(\delta_{1}, \delta_{2}) \in 
\{(0,1), (1,0), (1,1)\}$. 
However this is against (3.2). 

\par 
Next assume that $[\overline{\Delta} : \Delta] = 1$, that is, 
$\Delta$ is primitive in $\text{Pic} S_{3}$. Then, as before, there exists 
an element $h \in \text{Pic} S_{3}$ such that $\text{Pic} S_{3} = \Delta + \Bbb Z\cdot h$. 
Set $\Bbb Z \cdot H = \Delta^{\bot}$ in $\text{Pic} S_{3}$ and 
$n = [\text{Pic} S_{3} : \Delta \oplus \Bbb Z \cdot H]$. 
Then 
$n^{2} = ((\text{discr}\Delta) \cdot H^{2})/(\text{discr Pic} S_{3}) = 
40 H^{2}/3$ and (by replacing $h$ by $-h$ if necessary,) we can find 
integers $a_{i}$, $b_{j}$ such that 
$H = nh + \sum a_{i}C_{i} + \sum b_{j}E_{j}$ in $\text{Pic} S_{3}$, that is, 
$\frac{H}{n} = h + \sum \frac{a_{i}}{n}C_{i} + \sum \frac{b_{j}}{n}E_{j}$ 
in $\text{Pic} S_{3} \otimes \Bbb Q$. Using $H.C_{\alpha} = H.E_{\beta} = 0$ and 
the negative definiteness of $(C_{i}\cdot C_{\alpha})$ and 
$(E_{j}\cdot E_{\beta})$, we see that $(\frac{a_{i}}{n}, \frac{b_{j}}{n})$ 
is the unique solution of 
$$
(h + \sum x_{i}C_{i} + \sum y_{j}E_{j})\cdot C_{\alpha} = 0, 
(h + \sum x_{i}C_{i} + \sum y_{j}E_{j})\cdot E_{\beta} = 0, 
$$ 
($\alpha = 1, ..., 10; \beta = 1, ..., 9$). 

\par 
Since $\text{discr}(C_{i}\cdot C_{\alpha}) = 4$ and 
$\text{discr}(E_{j}\cdot E_{\beta}) = 10$, 
this implies that $\frac{a_{i}}{n} \in \frac{\Bbb Z}{4}$ and 
$\frac{b_{j}}{n} \in \frac{\Bbb Z}{10}$. 
Thus, 
$\frac{20H}{n} = 20h + \sum \frac{20a_{i}}{n}C_{i} 
+ \sum \frac{20b_{j}}{n}E_{j} \in \text{Pic} S_{3}$ 
and in particular, $n \vert20$. 
Combining this with $H^{2} = 3 n^{2}/40 \equiv 0 (\text{mod} 2)$, 
we find that $n = 20$ and $H^{2} = 30$. Thus, 
$H = 20h + \sum a_{i}C_{i} + \sum b_{j}E_{j}$, 
$a_{i} \in 5\cdot \Bbb Z$ and $b_{j} \in 2\cdot \Bbb Z$. 
Set $a_{i} = 5\alpha_{i}$ and $b_{j} = 2\beta_{j}$. 
Then, 
\roster 
\item 
$H = 20h + 5(\sum \alpha_{i}C_{i}) + 2(\sum \beta_{j}E_{j})$, and 
\item 
$400h^{2} = (20h)^{2} = H^{2} + 25(\sum \alpha_{i}C_{i})^{2} 
+ 4(\sum \beta_{j}E_{j})^{2}$, 
\item 
$0 = H \cdot (\sum \alpha_{i}C_{i}) = 20(h \cdot \sum \alpha_{i}C_{i}) 
+ 5(\sum \alpha_{i}C_{i})^{2}$, and 
\item 
$0 = H \cdot (\sum \beta_{j}E_{j}) = 20(h \cdot \sum \beta_{j}E_{j}) 
+ 2(\sum \beta_{j}E_{j})^{2}.$
\endroster
By (3) and (4), we see that 
$5(\sum \alpha_{i}C_{i})^{2}) \equiv 2(\sum \beta_{j}E_{j})^{2}) 
\equiv 0 (\text{mod} 20)$. 
Substituting this into (2), we get $H^{2} \equiv 0 (\text{mod} 20)$. 
However this contradicts the previous equality $H^{2} = 30$. 
Now we are done.
\qed 
\enddemo

\par
\head 
\S 4. Classification of extremal log Enriques surfaces
\endhead 

\par
In this section, we prove the Main Theorem (3). 
Throughout this section, we again work in the setting (2.1). 
By the Main Theorem (1), we know that $\Delta$ is now one of either 
$A_{19}$, $D_{19}$, $D_{16} \oplus A_{3}$, $D_{13} \oplus A_{6}$, 
$D_{7} \oplus A_{12}$, $D_{4} \oplus A_{15}$, or $D_{7} \oplus D_{12}$. \par 
In the case where $\Delta$ is either of type $A_{19}$ or of type $D_{19}$, 
the result follows from [OZ1, Theorems 1 and 2]. 
So we may consider the remaining cases: 

\par 
{\it Case} 1. $D_{16} \oplus A_{3}$, {\it Case} 2. $D_{13} \oplus A_{6}$, 
{\it Case} 3. $D_{7} \oplus A_{12}$, {\it Case} 4. $D_{4} \oplus A_{15}$,
and \par 
{\it Case} 5. $D_{7} \oplus D_{12}$. 

\par 
Since in each case $(S, g) \simeq (S_{3}, g_{3})$ ((2.1)), 
we identify these two in the sequel. 
Set $\Delta = (\cup_{i=1}^{3l+1} C_{i}) \coprod 
(\cup_{j=1}^{3m} E_{j})$ where in Cases (1) - (5),  
$\sum C_{i}$ is of type $D_{3l+1}$, 
$C_{3l+1}.C_{3l-1} = C_{3l}.C_{3l-1} = C_{i}.C_{i+1} = 1 (i \leq 3l-1)$, 
$\sum E_{j}$ is of type $A_{3m}$ and  
$E_{j}.E_{j+1} = 1 (j \geq 1)$, and in Case (5), 
$\sum C_{i}$ is of type $D_{7}$, $C_{7}.C_{5} = C_{6}.C_{5} = C_{j}.C_{j+1} = 1 (j \leq 4)$, $\sum E_{j}$ is of type $D_{12}$ and  
$E_{12}.E_{10} = E_{11}.E_{10} = E_{j}.E_{j+1} = 1 (j \leq 9)$.
We also denote by the same letter $\Delta$ the sublattice of 
$\text{Pic} S_{3}$ generated by the irreducible components of $\Delta$
and by $\overline{\Delta}$ its primitive closure in $\text{Pic} S_{3}$
as in Section 3.
Set 
$\Bbb Z \cdot H = \overline{\Delta}^{\bot} = \Delta^{\bot}$ in 
$\text{Pic} S_{3}$. Here we may take $H$ as the pull back of the ample generator 
of $\text{Pic} \overline{S}$. For convenience of notation, we sometimes set 
$G_{1} = C_{1}, ..., G_{3l+1} = C_{3l+1}, G_{3l+2} = 
E_{1}, ..., G_{19} = E_{3m}$. 

\par
\proclaim{Proposition (4.1)} In each case, $Z$ is unique up to isomorphisms 
if the following two conditions are satisfied: 
\roster
\item 
$H^{2}$ is determined only by the type of $Z$ and is 
independent of particular choice of $Z$.

\item 
The dual graph of the divisor  $\Delta$  determines,
uniquely up to graph isomorphisms,
rational numbers  $a_{ij}$, $b, b_{k}$ ($1 \leq i, j, k \leq 19$),
such that
$e_{j} := \sum_{i=1}^{19} a_{ij}G_{i}$ ($1 \leq j \leq 19$) form 
a  $\Bbb Z-$basis of $\overline{\Delta}$ and 
$e_{j}$'s, $e_{20} := bH + \sum_{k=1}^{19} b_{k}D_{k}$  form a 
$\Bbb Z-$basis of $\text{Pic} S_{3}$. 
\endroster
\endproclaim 

\par
\demo{Proof} Let  $Z$  be an extremal log Enriques surface
with  $\Delta, G_i, H$  as defined above or (2.1).
Let $Z(\alpha)$ be the extremal log Enriques surface
in (1.7) of the same type as that of  $Z$.
As for  $Z$, we can define similarly
$\Delta(\alpha)$, $G_{i}(\alpha)$, $H(\alpha)$, etc.
Then, by the conditions (1) and (2), there exists an isometry 
$\psi : \text{Pic} S_{3} \rightarrow \text{Pic} S_{3}$ such that 
$\psi(D_{i}) = D_{i}(\alpha)$, $\psi(H) = H(\alpha)$ and that 
$\psi$ preserves ample classes. The last condition follows from 
the fact, which is derived from Kleiman's criterion on ampleness, 
that there are sufficiently small positive 
numbers  $\gamma_k$  such that both  $H - \sum \gamma_{k}D_{k}$ and 
$H(\alpha) - \sum \gamma_{k}D_{k}(\alpha)$  are ample divisors. 
Then by [V, page 13], $\psi$ 
extends to an effective Hodge isometry $\overline{\psi}$ of 
$H^{2}(S_{3}, \Bbb Z)$.  Now we may apply the Torelli Theorem 
for K3 surfaces to get an automorphism $\varphi$ of 
$S_{3}$ such that $\varphi^{*} = \overline{\psi}$. By construction, 
$\varphi$ maps the exceptional divisor $\Delta$ to $\Delta(\alpha)$. 
Combining this with the result  $g_{3}\circ \varphi = \varphi \circ g_{3}$
in (1.2), we 
see that $\varphi$ is an equivariant isomorphism between the 
triplets $(S_{3}, g_{3}, \Delta)$ and $(S_{3}, g_{3}, \Delta(\alpha))$. 
Thus $\varphi$ descends to an isomorphism $Z \rightarrow Z(\alpha)$.  
\qed
\enddemo 

\par
Now we may check the conditions (1) and (2) for each case. \par
{\it Case 1, the case where} $\Delta$ {\it is of type} $D_{16} \oplus A_{3}$ \par 
\proclaim{Claim (4.2)} 
\roster 
\item 
$[\overline{\Delta} : \Delta] = 2$. 
\item 
Up to $Aut_{graph}(D_{16})$,
$e_{1} := \frac{1}{2}(C_{1} + C_{3} + C_{5} + \cdots + C_{13} + C_{16})$, 
$e_i = C_{17-i}$ ($2 \le i \le 16$),
$e_{17} := E_{1}$, $e_{18} := E_{2}$, $e_{19} := E_{3}$ 
form a  $\Bbb Z-$basis of $\overline{\Delta}$.
\endroster
\endproclaim 

\par
\demo{Proof} Since $[\overline{\Delta} : \Delta]^2 
= (\text{discr} \Delta)/(\text{discr} \overline{\Delta}) 
= 16/(\text{discr} \overline{\Delta})$,  
$[\overline{\Delta} : \Delta]$  is either  1, 2 or 4.

\par
If  $[\overline{\Delta} : \Delta] = 1$  or 4, we will 
get a contradiction as in Section 3, noting that we also
have  $\text{discr} (C_{i} \cdot C_{\alpha}) = \text{discr}
(E_{j} \cdot E_{\beta}) = 4$  here.
This proves (1) of (4.2).

\par
The fact that  $[\overline{\Delta} : \Delta] = 2$
and the argument in (3.1) imply that 
$\overline{\Delta} - \Delta$  contains, after interchanging
$C_1$  and  $C_2$  by the non-trivial element in
$\text{\rm Aut}_{\text{\rm graph}}(D_{16})$  if necessary,
either
$L = \frac{1}{2}(C_{1} + C_{3} + \cdots + C_{13} + C_{16})$,  or
$L = \frac{1}{2}(E_{1} + E_{3})$,  or 
$L = \frac{1}{2}(C_{1} + C_{3} + \cdots + C_{13} + C_{16}) + 
\frac{1}{2}(E_{1} + E_{3})$. 
Combining this with Nikulin's result (cf. Lemma (3.2)), we get 
$L = \frac{1}{2}(C_{1} + C_{3} + \cdots + C_{13} +
C_{16}) + \frac{1}{2}(E_{1} + E_{3})$. This implies the assertion (2). 
\qed 
\enddemo 

\par
\proclaim{Claim (4.3)}
\roster
\item 
$H^{2} = 12$. 
\item 
Up to $Aut_{graph}(\Delta)$, 
$e_{1}$, $e_{2}$, ..., $e_{19}$ and $\frac{1}{4}(H - E_{1} - 2E_{2} - 3E_{3})$ 
form a  $\Bbb Z-$basis of $\text{Pic} S_{3}$. 
\endroster
\endproclaim 
Set $n = [\text{Pic} S_{3} : \overline{\Delta} \oplus \Bbb Z \cdot H]$ and 
$\text{Pic} S_{3} = \overline{\Delta} + \Bbb Z \cdot h$. Then, we have
$n^{2} = (\text{discr} \overline{\Delta} \cdot H^{2})/(\text{discr Pic} S_{3}) 
= 4H^{2}/3$, that is, $H^{2} = 3 n^2/4$ and replacing 
$h$ by $-h$ if necessary we can write 
$H = nh + \sum \alpha_{k}e_{k}$ (for some integers $\alpha_{k}$). 
Using $\text{discr}(e_{i} \cdot e_{\alpha})_{i, \alpha = 1, ..., 16} = 1$ and 
$\text{discr}(e_{i} \cdot e_{\alpha})_{i, \alpha = 17, ..., 19} = 4$, 
we see that 
$\frac{\alpha_{i}}{n} \in \Bbb Z$ for $i = 1, ..., 16$ and 
$\frac{\alpha_{i}}{n} \in \frac{\Bbb Z}{4}$ for $i = 17, ..., 19$. 
In particular, 
$\frac{4H}{n} = 4h + \sum \frac{4\alpha_{k}}{n}e_{k} \in \text{Pic} S_{3}$. 
Thus $n \vert 4$. Combining this with 
$H^{2} = 3 n^{2}/4 \equiv 0 (\text{mod} 2)$, we find that 
$n = 4$, $H^{2} = 12$.  Thus, there exist integers 
$c_{i}$ ($i = 1, ..., 16$), $d_{j}$ ($j = 17, ..., 19$) such that 
$H = 4h + \sum_{i=1}^{16} 4c_{i}e_{i} + \sum_{j=17}^{19} d_{j}e_{j}$. 
Replacing $h$ by $h - \sum m_{k}e_{k}$ ($m_{k} \in \Bbb Z$), and 
using $e_{j+16} = E_{j}$ if $j = 1, ..., 3$, we can adjust $h$ like 
$h = \frac{1}{4}H - \frac{1}{4} \sum_{j = 1}^{3} a_{j}E_{j}$ for some integers $a_{j} \in \{0, 1, 2, 3\}$. We shall determine $a_{j}$ up to 
$\text{\rm Aut}_{\text{\rm graph}}(A_{3})$. 
Using $h \cdot E_{j} \in \Bbb Z$, we get 
$-2a_{1} + a_{2} \equiv 0 (\text{mod} 4)$, 
$a_{1} - 2a_{2} + a_{3} \equiv 0 (\text{mod} 4)$, and 
$a_{2} - 2a_{3} \equiv 0 (\text{mod} 4)$. Thus, up to 
$\text{\rm Aut}_{\text{\rm graph}}(A_{3})$, 
we have either 
(1) $a_{1} = a_{2} = a_{3} = 0$, or (2) $a_{1} = a_{3} = 2$ and $a_{2} = 0$ 
or (3) $a_{1} = 1$, $a_{2} = 2$ and $a_{3} = 3$. 

\par 
In case (1), we calculate 
$h^{2} = (\frac{1}{4}H)^{2} = \frac{3}{4} \not\in \Bbb Z$, 
a contradiction. Also, in case (2), we calculate 
$h^{2} = (\frac{1}{4}H)^{2} + \frac{1}{4}(E_{1} + E_{3})^{2}  
= \frac{3}{4} - 1 \not\in \Bbb Z$, a contradiction. 
Thus, the only possible values of $a_{j}$ are $a_{1} = 1$, $a_{2} = 2$ and $a_{3} = 3$. 
Since we already know the existence of such $a_{j}$, this gives the 
assertion (2). 
\qed 

\par
\remark{Remark} It does not seem easy at least for the authors to 
find directly for which $a_{j}$,
$\frac{1}{4}H - \frac{1}{4} \sum_{j = 1}^{3} a_{j}E_{j}$ is really 
in $\text{Pic} S_{3}$. This is the reason why we argued as above.
\endremark 

\par \vskip 5pt
{\it Case 2, the case where} $\Delta$ {\it is of type} $D_{13} 
\oplus A_{6}$ 

\par 
\proclaim{Claim (4.4)} The sublattice  $\Delta$  is primitive 
in $\text{Pic} S_{3}$, i.e., $\Delta$  is equal to  its
primitive closure  $\overline{\Delta}$  in  $Pic S_3$.
\endproclaim 
\demo{Proof} Since $\text{discr} \Delta = 4 \cdot 7$, if (4.4) is false
we have  $[\overline{\Delta} : \Delta] = 2$. Then we will reach
a contradiction to (3.2) as in the proof of (3.1).
\qed
\enddemo

\proclaim{Claim (4.5)} 
\roster 
\item 
$H^{2} = 84$ and 
\item 
Up to $Aut_{graph}(\Delta)$, 
$C_{1}, ..., C_{12}, E_{1}, ..., E_{6}$ and 
$\frac{1}{28}H - 
\frac{1}{4}(2C_{1} + 2C_{3} + \cdots + 2C_{9} + 2C_{11} + C_{12} + 3C_{13}) - \frac{1}{7}(\sum_{j=1}^{6} jE_{j})$ 
are $\Bbb Z-$basis of $\text{Pic} S_{3}$. 
\endroster
\endproclaim 

\par
\demo{Proof} Set $\text{Pic} S_{3} = \Delta + \Bbb Z \cdot h$ and
$n = [\text{Pic} S_{3} : \Delta \oplus \Bbb Z\cdot H]$. 
Then by the same argument as before, we see that 

\par 
$n^{2} = (\text{dicsr} \Delta \cdot H^{2})/(\text{discr Pic} S_{3}) = 
28H^{2}/3$, 

\par 
$H = nh + \sum_{i} \alpha_{i}C_{i} + \sum_{j} \beta_{j}E_{j}$, and \par 
$\frac{\alpha_{i}}{n} \in \frac{\Bbb Z}{4}$, 
$\frac{\beta_{j}}{n} \in \frac{\Bbb Z}{7}$. 
Thus $\frac{28}{n}H = 28h + 7(\sum \frac{4\alpha_{i}}{n}C_{i}) 
+ 4(\sum \frac{7\beta_{j}}{n}E_{j}) \in \text{Pic} S_{3}$ and then 
$n \vert 28$. Combining this with 
$H^{2} = \frac{3}{28}n^{2} \equiv 0 (\text{mod} 2)$, we get 
$n = 28$ and $H^{2} = 84$. Thus, replacing $h$ by 
$h - \sum_{i} m_{i}C_{i} - \sum_{j} n_{j}E_{j}$ 
($m_{i}, n_{j} \in \Bbb Z$), we may adjust $h$ such as 
$h = \frac{1}{28}H - \frac{1}{4}(\sum_{i} a_{i}C_{i}) - 
\frac{1}{7}(\sum_{j} b_{j}E_{j})$ where  
$a_{i} \in \{0,1,2,3\}$ and $b_{j} \in \{0,1, ..., 6\}$. 
We determine $a_{i}$ and $b_{j}$ up to $Aut_{graph}(\Delta)$. 
By $Aut_{graph}(\Delta)$, we may assume $a_{12} \leq a_{13}$ 
and $b_{1} \leq b_{6}$. Using $h \cdot C_{k} \in \Bbb Z$, we can readily 
see that $a_{i} \equiv ia_{i} (\text{mod} 4)$ for $1 \leq i \leq 11$, 
$a_{10} + a_{12} + a_{13} - 2a_{11} \equiv 0 (\text{mod} 4)$, 
$-2a_{12} + a_{11} \equiv 0 (\text{mod} 4)$, and 
$-2a_{13} + a_{11} \equiv 0 (\text{mod} 4)$. 
These formulas imply that $a_{i}$ are either, 

\roster
\item 
$a_{1} = a_{2} = \cdots = a_{13} = 0$, 
\item 
$a_{1} = a_{2} = \cdots = a_{11} = 0$, $a_{12} = a_{13} = 2$, or 
\item 
$a_{1} = a_{3} = \cdots = a_{9} = a_{11} = 2$, 
$a_{2} = a_{4} = \cdots = a_{10} = 0$, $a_{12} = 1$ and $a_{13} = 3$. 
\endroster 

\par
Using $h \cdot E_{l} \in \Bbb Z$, we have 
$b_{j} \equiv jb_{1} (\text{mod} 7)$ and $b_{1}$ is either $0, 1, 2, 3$.
The assertion here follows from our assumption that
$b_{1} \leq b_{6}$.  Thus, according to (1), (2), (3), 
$h^{2} \equiv \frac{3}{28} + \frac{b_1^2}{7}$ ($\text{mod} \Bbb Z$), 
$\equiv \frac{3}{28} - 1 + \frac{b_1^2}{7}$ ($\text{mod} \Bbb Z$), and 
$\equiv \frac{3}{28} - \frac{13}{4} + \frac{b_1^2}{7}$ ($\text{mod} \Bbb Z$). 
Thus, by $h^{2} \equiv 0$ ($\text{mod} 2\cdot\Bbb Z$), we see that 
$b_{1} = 1$ and that $a_{i}$ satisfy (3). This proves (4.5). 
\qed 
\enddemo 

\par \vskip 5pt
{\it Case 3, the case where} $\Delta$ {\it is of type} $D_{7} 
\oplus A_{12}$ 

\par 
By the same argument as in case 2, we get the following two claims, which 
guarantee the conditions (1) and (2) in (4.1).

\par
\proclaim{Claim (4.6)} $\Delta$ is primitive in $\text{Pic} S_{3}$.
\endproclaim 
\proclaim{Claim (4.7)} 
\roster 
\item 
$H^{2} = 156$ and 
\item 
Up to $Aut_{graph}(\Delta)$, 
$C_{1}, ..., C_{7}, E_{1}, ..., E_{12}$ and 
$\frac{1}{52}H - 
\frac{1}{4}(2C_{1} + 2C_{3} + 2C_{5} + C_{6} + 3C_{7}) - \frac{1}{7}(\sum_{j=1}^{12} \overline{2j}E_{j})$ 
are $\Bbb Z-$basis of $\text{Pic} S_{3}$, where we denote by $\overline{2j}$ the integer determined by $2j \equiv \overline{2j} (\text{mod} 13)$ and 
$0 \leq \overline{2j} \leq 12$.
\endroster
\endproclaim  

\par \vskip 5pt
{\it Case 4, the case where} $\Delta$ {\it is of type} $D_{4} \oplus A_{15}$

\par 
This is the hardest case. 

\par
\proclaim{Claim (4.8)}
\roster
\item
$[\overline{\Delta} : \Delta] = 4$ and 
$\overline{\Delta}/\Delta \simeq \Bbb Z/4$.
\item 
$e_{1} := C_{1}$, $e_{2} := C_{2}$, $e_{3} := C_{3}$, $e_{4} := C_{4}$, 
$e_{5} := E_{2}$, $e_{6} := E_{3}$, $\cdots$ $e_{17} := E_{14}$, 
$e_{18} := E_{15}$ and 
$e_{19} := \frac{1}{2}(C_{1} + C_{2}) + 
\frac{1}{4}(E_{1} + 2E_{2} - E_{3} + E_{5} + 2E_{6} - E_{7} + E_{9} 
+ 2E_{10} - E_{11} + E_{13} + 2E_{14} - E_{15})$ are $\Bbb Z-$basis 
of $\overline{\Delta}$. 
\endroster 
\endproclaim 

\par
\demo{Proof of (1)} 
Since $\text{discr} \Delta = 4 \cdot 16 = 64$, we have either 
(i) $[\overline{\Delta} : \Delta] = 8$, 
(ii) $[\overline{\Delta} : \Delta] = 1$,
(iii) $[\overline{\Delta} : \Delta] = 2$, 
(iv) $[\overline{\Delta} : \Delta] = 4$ and 
$\overline{\Delta}/\Delta \simeq (\Bbb Z/2)^{\oplus 2}$, or 
(v) $[\overline{\Delta} : \Delta] = 4$ and 
$\overline{\Delta}/\Delta \simeq \Bbb Z/4$. 

\par 
We elminate the cases (i) - (iv) by arguing by contradiction. \par 
{\it Case (i).}  In this case $\text{discr} \overline{\Delta} = 1$. Then, 
$H^{2} = \text{discr Pic} S_{3} = 3 \not\equiv 0 (\text{mod} 2)$, a contradiction. 

\par 
{\it Case (ii).} We have $\overline{\Delta} = \Delta$. 
Set $n = [\text{Pic} S_{3} : \Delta \oplus \Bbb Z \cdot H]$. Then, 
$n^{2} = \frac{64}{3}H^{2}$ and $H = nh + \sum_{i} a_{i}C_{i} 
+ \sum_{j} b_{j}E_{j}$ for some integers $a_{i}, b_{j}$. 
Since $\text{discr} (C_{i} \cdot C_{\alpha}) = 4$ and 
$\text{discr} (E_{j} \cdot E_{\beta}) = 16$, we see that 
$\frac{4a_{i}}{n}, \frac{16b_{j}}{n} \in \Bbb Z$ and that 
$\frac{16}{n}H = 16h + \sum \frac{16a_{i}}{n}C_{i} + \sum \frac{16b_{j}}{n}E_{j} 
\in \text{Pic} S_{3}$. Thus $n \vert 16$. 
Combining this with $H^{2} = \frac{3}{64}n^{2} \equiv 0 (\text{mod} 2)$, we get 
$n =16$ and $H^{2} = 12$. Then 
$H = 16h + \sum_{i} 4\alpha_{i}C_{i} + \sum_{j} b_{j}E_{j}$, 
where $\alpha_{i} = \frac{a_{i}}{4} (\in \Bbb Z)$. Using this formula, 
we calculate 
$16^{2}h^{2} = H^{2} + 16(\sum_{i} \alpha_{i}C_{i})^{2} 
+ (\sum_{j} b_{j}E_{j})^{2}$. 
On the other hand, since 
$0 = H \cdot (\sum_{j} b_{j}E_{j}) = 16h \cdot (\sum_{j}b_{j}E_{j}) 
+ (\sum_{j} b_{j}E_{j})^{2}$, we find 
$(\sum_{j} b_{j}E_{j})^{2} \equiv 0 (\text{mod} 16)$. 
Then $12 = H^{2} \equiv 0 (\text{mod} 16)$, a contradiction. 

\par 
{\it Case (iii).}  In this case, there exist integers 
$\alpha_{i}, \beta_{j} \in \{0,1\}$ such that 
$L:= \frac{1}{2}(\sum_{i=1}^4 \alpha_{i}C_{i} 
+ \sum_{j = 1}^{15} \beta_{j}E_{j}) \in \overline{\Delta} - \Delta$. 
Since $L.C_{i} \in \Bbb Z$ and $L.E_{j} \in \Bbb Z$, 
we readily find that 
$L = \frac{1}{2}(E_{1} + E_{3} + \cdots + E_{13} + E_{15})$ 
and that $C_{1}, ..., C_{4}, G_{1} := E_{1}, G_{2} := E_{2}, ..., 
G_{14} := E_{14}$, and $G_{15} := L$ are $\Bbb Z-$basis of 
$\overline{\Delta}$. Set 
$n = [\text{Pic} S_{3} : \overline{\Delta} \oplus \Bbb Z \cdot H]$ 
and $\text{Pic} S_{3} = \overline{\Delta} + \Bbb Z \cdot h$. 
Then $n^{2} = \frac{16}{3}H^{2}$ and 
$H = nh + \sum_{i} a_{i}C_{i} + \sum_{j} b_{j}G_{j}$ for some integers
$a_{i}, b_{j}$. Since 
$\text{discr} (C_{i} \cdot C_{\alpha}) = \text{discr} (G_{j} \cdot G_{\beta}) = 4$, 
we have $\frac{a_{i}}{n}, \frac{b_{j}}{n} \in \frac{\Bbb Z}{4}$, that is, 
$\frac{4H}{n} = 4h + \sum_{i} \frac{4a_{i}}{n}C_{i} 
+ \sum_{j} \frac{4b_{j}}{n}G_{j} \in \text{Pic} S_{3}$. 
Thus $n \vert 4$. Then $H^{2} = \frac{3}{16}n^{2} \not\equiv 0 (\text{mod} 2)$, 
a contradiction. 

\par 
{\it Case (iv).} In this case there should exist at least two 
$L_{1}, L_{2} \in \overline{\Delta} - \Delta$. 
However, the same argument as in case (3) shows that such $L_{i}$ is unique, namely 
$\frac{1}{2}(E_{1} + E_{3} + \cdots + E_{13} + E_{15})$,
a contradiction. \par 
Now the assertion (1) is proved. 
\qed 
\enddemo 

\par
\demo{Proof of (2)} By (1), there exist 
subsets $I \subset \{1, 2, 3, 4\}$, $J \subset \{1, 2, ...,15\}$ 
and integers $\alpha_{i} \in \{1, 2, -1\}$,  
$\beta_j \in \{1, 2, -1\}$ such that  
$N := \frac{1}{4}(\sum_{i \in I} \alpha_{i}C_{i} + 
\sum_{j \in J} \beta_{j}E_{j}) \in \overline{\Delta} - \Delta$ 
and that $2N \not\in \Delta$. We determine $N$. 
Set $I' := \{i \in I \vert \alpha_{i} \not= 2 \}$,  
$J' := \{j \in J \vert \beta_{j} \not= 2 \}$ and 
$N' := \sum_{i \in I'} \vert \alpha_{i} \vert C_{i} 
+ \sum_{j \in J'} \vert \beta_{j} \vert E_{j}$.
Then
$0 \equiv 4N \equiv N' (\text{mod} 2\cdot\text{Pic} S_{3})$. 
On the other hand, since $2N \not\in \Delta$, $I' \not= \phi$ or 
$J' \not= \phi$. 
Using $N' \cdot C_{i} \equiv N' \cdot E_{j} \equiv 0 (\text{mod} 2)$ and (3.2), 
we find that $I' = \phi$ and $J' = \{1, 3, 5, ..., 13, 15 \}$. 
Replacing $N$ by $-N$ if necessary, we may assume that $\beta_{1} = 1$. 
Set $M := E_{1} + 2E_{2} - E_{3} + E_{5} + 2E_{6} - E_{7} + E_{9} 
+ 2E_{10} - E_{11} + E_{13} + 2E_{14} - E_{15}$. 
Then using $N \cdot C_{i} \in \Bbb Z$ and $N \cdot E_{j} \in \Bbb Z$, 
we readily see that (up to $Aut_{graph}(\Delta)$,) $N$ is either 
(1) $\frac{1}{2}(C_{1} + C_{2}) + \frac{1}{4}M$ or (2) $\frac{1}{4}M$. 
However, in case (2), 
$N^{2} = \frac{1}{16}M^{2} = -3 \not\equiv 0 (\text{mod} 2)$, a contradiction. 
Thus $N = \frac{1}{2}(C_{1} + C_{2}) + \frac{1}{4}M = e_{19}$. 
This implies the assertion (2). 
\qed 
\enddemo 

\par
\proclaim{Claim (4.9)} 
\roster 
\item 
$H^{2} = 12$. 
\item 
(up to $Aut_{graph}(\Delta)$,)
$e_{1}$, $e_{2}$, ..., $e_{19}$ and 
$e_{20} := \frac{1}{4}(H - e_{1} - 3e_{2} - 2e_{4}  
- 2e_6 - 2e_7 - 2e_8 - 2e_9
- 2e_{14} - 2e_{15} - 2e_{16} 
- 2e_{17} - 2e_{19})$ are $\Bbb Z-$basis of $\text{Pic} S_{3}$. 
\endroster 
\endproclaim 

\par
\demo{Proof} Set $n = [\text{Pic} S_{3} : \overline{\Delta} \oplus \Bbb Z \cdot H]$ and $\text{Pic} S_{3} = \overline{\Delta} + \Bbb Z \cdot h$.
Then using the same argument as in case 1 based on 
$\text{discr} \overline{\Delta} = 4$, we get $n = 4$, $H^{2} = 12$ and 
find integers $a_{k} \in \{0, 1, 2, 3\}$ ($1 \leq k \leq 19$) such that 
$e_{1}, ..., e_{19}$ and 
$e := \frac{1}{4}(H - \sum_{k} a_{k}e_{k})$ are $\Bbb Z-$ basis of 
$\text{Pic} S_{3}$. We determine $a_{i}$ up to $Aut_{graph}(\Delta)$.  Since 
$(\sum_{k} a_{k}e_{k}) \cdot e_{l} = -4e \cdot e_{l} \equiv 0 (\text{mod} 4)$, 
we see that $a_{k}$ are either 

\roster
\item 
$a_{19} = a_{18} = \cdots a_{5} = a_{4} = \cdots a_{1} = 0$, 
\item 
$a_{19} = a_{18} = \cdots a_{5} = 0$, $a_{4} = \cdots = a_{1} = 2$, or 
\item 
$a_{19} = 2$, $a_{18} = 0$, $a_{17} = \cdots a_{14} = 2$, 
$a_{13} = \cdots a_{10} = 0$, $a_{9} = \cdots = a_{6} = 2$, $a_{5} = 0$, 
$a_{4} = 2$, $a_{3} = 0$, $a_{2} = 3$, and $a_{1} = 1$ (up to 
$Aut_{graph}(D_{4})$). 
\endroster 

\par
However, in cases of (1) and (2), we see that $e^{2} \not\in \Bbb Z$, 
a contradiction. Thus the case (3) occurs, that is, $e = e_{20}$. 
\qed 
\enddemo 

\par \vskip 5pt
{\it Case 5, the case where} $\Delta$ {\it is of type} $D_{7} 
\oplus D_{12}$ 

\par 
The verification is also quite similar. 
We only indicate Claims needed to check the conditions (1) and (2) in (4.1). 
\proclaim{Claim (4.10)}
\roster
\item 
$[\overline{\Delta} : \Delta] = 2$. 
\item 
(Up to $Aut_{graph}(\Delta)$) 
$e_{1} := C_{1}, e_{2} := C_{2}, ..., e_{6} := C_{6}$, 
$e_{7} := \frac{1}{2}(C_{6} + C_{7} + E_{1} + E_{3} + \cdots + E_{9} + E_{11})$,  
$e_{8} := E_{1}$, $e_{9} := E_{2}$, ..., $e_{19} := E_{12}$ are 
$\Bbb Z-$basis of $\overline{\Delta}$. 
\endroster
\endproclaim

\par
\proclaim{Claim (4.11)} 
\roster
\item 
$H^{2} = 12$. 
\item 
Up to an element of $Aut_{graph}(\Delta)$ which keeps $e_{1}, ..., e_{19}$ invariant,  
$e_{1}, ..., e_{19}$ and 
$e_{20} := \frac{1}{4}(H - 2e_{1} - 2e_{2} - \cdots - 2e_{7} - 
e_{8} - e_{10} - \cdots - e_{16} + e_{18} - e_{19})$ are 
$\Bbb Z-$basis of $\text{Pic} S_{3}$.
\endroster
\endproclaim 

\par
Now we have completed the proof of the Main Theorem (3). Q.E.D. 

\par 
 
\enddocument